\documentclass[12pt]{article}

\usepackage{amsmath,framed}
\usepackage{amssymb}
\usepackage{braket} %Set notation
\usepackage{enumerate} % roman enumeration
\usepackage{amsthm}
\usepackage{graphicx}
\usepackage{psfrag}
\usepackage[T1]{fontenc}

\newtheorem{theorem}{Theorem}[section]
\newtheorem{proposition}[theorem]{Proposition}
\newtheorem{lemma}[theorem]{Lemma}
\newtheorem{corollary}[theorem]{Corollary}

\numberwithin{equation}{section}

\theoremstyle{definition}
\newtheorem{definition}[theorem]{Definition}
\newtheorem{remark}[theorem]{Remark}

\newtheorem{question}[theorem]{Question}

\newcommand{\IIi}{${\rm{II}}_1\ $}

\newcommand{\nm}[1]{\left\|{#1}\right\|}
\newcommand{\ip}[2]{\left<#1,#2\right>}
\newcommand{\F}{\mathbb F}
\newcommand{\LF}{\mathcal L\F}
\newcommand{\vnotimes}{\overline{\otimes}}

\title{\textsc{The radial masa in a free group factor is maximal injective}}
\author{Jan Cameron\\\normalsize{\texttt{jacameron@vassar.edu}}\and Junsheng Fang\thanks{Partially supported by the fundamental research funds for the central universities of China.}\\\normalsize{\texttt{junshengfang@gmail.com.}}\and Mohan Ravichandran\\\normalsize{\texttt{mohanr@sabanciuniv.edu}}\and Stuart White\\\normalsize{\texttt{stuart.white@glasgow.ac.uk}}}
\begin{document}
\maketitle

\begin{abstract}
The radial (or Laplacian) masa in a free group factor is the abelian
von Neumann algebra generated by the sum of the generators (of
the free group) and their inverses. The main result of this paper is that the radial masa is a maximal injective von Neumann subalgebra of a free group factor. We also investigate tensor products of maximal injective algebras. Given two inclusions $B_i\subset M_i$ of type $\mathrm{I}$ von Neumann algebras in finite von Neumann algebras such that each $B_i$ is maximal injective in $M_i$, we show that the tensor product $B_1\ \vnotimes\ B_2$ is maximal injective in $M_1\ \vnotimes\ M_2$ provided at least one of the inclusions satisfies the asymptotic orthogonality property we establish for the radial masa. In particular it follows that finite tensor products of generator and radial masas will be maximal injective in the corresponding tensor product of free group factors.
\end{abstract}

\section{Introduction}

In \cite{Popa.MaxInjective}, Popa showed that the von Neumann algebra generated by a single generator $a_1$ of a free group $\F_K=\langle a_1,\dots,a_K\rangle$ is a maximal injective von Neumann subalgebra of the factor $\LF_K$ associated to $\F_K$. This resolved a problem of Kadison by demonstrating that self-adjoint elements in a $\mathrm{II}_1$ factor need not contained in a hyperfinite subfactor. The subalgebra generated by $a_1$ is known as a generator masa (maximal abelian subalgebra) in $\LF_K$. There is another naturally occurring masa in the free group factor $\LF_K$, the von Neumann subalgebra generated by $\sum_{i=1}^K(a_i+a_i^{-1})$. This was shown to be maximal abelian by Pytlik in \cite{Pytlik.RadialMasa} and is known as the radial masa in $\LF_K$. The main objective of this paper is to show that this radial masa gives another example of an abelian maximal injective von Neumann subalgebra of $\LF_K$.

The radial masa shares many properties with generator masas. Dixmier defined a masa $A$ in a \IIi factor $M$ to be \emph{singular} if every unitary $u\in M$ with $uAu^*=A$ lies in $A$, and showed in \cite{Dixmier.Masa} that a generator masa in a free group factor is singular. Singularity of the radial masa was established by R\u{a}dulescu by means of an intricate calculation in \cite{Radulescu.LapMasa}, which provides a central ingredient in this paper. An alternative combinatorial proof was given by Sinclair and Smith in \cite{Sinclair.LapMasa}. To establish singularity of the radial masa, R\u{a}dulescu computed its Puk\'anszky invariant, historically the most sucessful invariant for singular masas (see \cite{Sinclair.MasaBook} for a discussion of this invariant). The Puk\'anszky invariant of a generator masa is easily computed using the group-subgroup methods of \cite{Sinclair.Puk} and both the radial and generator masas have Puk\'anszky invariant $\{\infty\}$. It is natural to ask whether these masas are conjugate in $\LF_K$. This seems unlikely, but we have been unable to find a proof. As noted in Proposition \ref{Inner} below, a straightforward deduction from \cite{Popa.Orth} shows that the radial masa is not inner conjugate to a generator masa in $\LF_K$. 
\begin{question}
Does there exist an automorphism of $\LF_K$ which maps the radial masa onto a generator masa?
\end{question}

The critical ingredient in Popa's proof of the maximal injectivity of a generator masa $B$ in a free group factor $\LF_K$ is an asymptotic orthogonality property at the Hilbert space level. Lemma 2.1 of \cite{Popa.MaxInjective} shows that if $x_1,x_2$ are elements of an ultrapower $(\LF_K)^\omega$ with $E_{B^\omega}(x_1)=E_{B^\omega}(x_2)=0$ and $y_1,y_2$ are elements of $\LF_K$ with $E_B(y_1)=E_B(y_2)=0$, then $x_1y_1$ is orthogonal to $y_2x_2$ in the Hilbert space $L^2((\LF_K)^\omega)$. To do this, Popa first shows that those elements $x\in\LF_K$ which approximately commute with the generator $a_1$ and are orthogonal to $B$ must be essentially supported on the collection of words which begin and end with large powers of $a_1$ or $a_1^{-1}$.  The asymptotic orthogonality follows when $y_1,y_2$ are group elements in $\F_K\setminus\{a_1^n:n\in\mathbb Z\}$ and a density argument then gives the result for general $y_1$ and $y_2$ (see also \cite[Section 14.2]{Sinclair.MasaBook}).

Most of this paper is taken up with establishing an asymptotic orthogonality condition for the radial masa (Theorem \ref{Radial.AOP}). The natural orthonormal basis of $\ell^2(\F_K)$ does not behave well in calculations involving the radial masa, however to prove that the radial masa is singular R\u{a}dulescu introduced in \cite{Radulescu.LapMasa} a collection of vectors which form a basis (albeit not quite an orthonormal basis) for $\ell^2(\F_K)$ which can be thought of as spanning double cosets coming from the radial masa. In Section \ref{Prelim}, we set up notation for discussing the radial masa and describe the properties of R\u{a}dulescu's basis.  In Section \ref{Asymptotics}, we establish the essential location of the support of an element $x$ in $\LF_K$ which is both orthogonal to the radial masa and approximately commutes with it. We do this in Lemma \ref{Asymptotics.Main}, the proof of which is contained in the lemmas preceding it.   Section \ref{Count} contains the second step of \cite[Lemma 2.1]{Popa.MaxInjective} for the radial masa. In this section, we use techniques from \cite{Sinclair.LapMasa} to count certain words in $\F_K$ which we use in Section \ref{Radial} to establish the asymptotic orthogonality result from the results of Section \ref{Asymptotics}.  All of the calculations described above are performed in $\LF_K\ \vnotimes\ N$, where $N$ is an arbitrary finite von Neumann algebra. The only additional difficulties this introduces are notational, and it enables us to show that certain tensor products involving the radial masa are also maximal injective, as described below.

In \cite[4.5 (1)]{Popa.MaxInjective}, Popa asked how maximal injective von Neumann algebras behave under tensor products. That is, if $B_i\subset M_i$ ($i=1,2$) are two inclusions of von Neumann algebars with $B_i$ a maximal injective von Neumann subalgebra of $M_i$, must $B_1\ \vnotimes\ B_2$ be a maximal injective von Neumann subalgebra of $M_1\ \vnotimes\ M_2$?  Various authors have subsequently worked on this question. The first progress was made by Ge and Kadison, who in \cite{Ge.Kadison} gave a positive answer when $B_1=M_1$ is an injective factor. This result was subsequently improved by Str\u{a}til\u{a} and Zsid\'o who removed the factor assumption to show (\cite[Theorem 6.7]{Stratila.Commutation}) that if $M_1$ is an injective von Neumann algebra and $M_2$ is a von Neumann algebra with a separable predual, then $M_1\ \vnotimes \ B_2$ is a maximal injective von Neumann subalgebra of $M_1\ \vnotimes\ M_2$ whenever $B_2$ is a maximal injective von Neumann subalgebra of $M_2$. In \cite{Shen.MaxInjective}, Shen investigated the tensor product 
of copies of a generator masa inside the a free group factor, answering Popa's question positively in this case. He was also able to obtain a positive result for an infinite tensor product of generator masas and so obtain the first example of an abelian subalgebra which is maximal injective in  a McDuff \IIi factor. In \cite{Shen.MaxInjective}, it was necessary to develop an additional orthogonality condition styled upon \cite[Lemma 2.1]{Popa.MaxInjective} for tensor products of generator masas.  Recent progress was made by the second author, who in \cite{Fang.MaxInjective} gave a positive answer to Popa's question in the case that $Z(B_1)$ is atomic and $M_1$ has a separable predual.  

In Section \ref{MaxInjective}, we examine how asymptotic orthogonality properties based on \cite[Lemma 2.1]{Popa.MaxInjective} imply maximal injectivity. Using Popa's intertwining lemma and ingredients from \cite{Fang.MaxInjective,Shen.MaxInjective}, we are able to give a technical improvement of the argument which deduces the maximal injectivity of the generator masa from the asymptotic orthogonality condition. This argument doesn't require any further assumptions on the masa beyond singularity and is also applicable to the tensor product of maximal injective algebras, provided we can establish our  asymptotic orthogonality condition after tensoring by a further finite von Neumann algebra (see Definition \ref{AOP} below).  In particular, it follows that if $A\subset M$ is a singular masa of a \IIi factor with the asymptotic orthogonality property, such as a generator masa or the radial masa in a free group factor, and $B$ is a type $\mathrm{I}$ maximal injective von Neumann subalgebra of a \IIi factor $N$ with separable predual, then $A\ \vnotimes\ B$ is maximal injective in $M\ \vnotimes\ N$ (Theorem \ref{MaxInjective.Final}, taking $A$ to be a masa). This result strengthens Shen's finite tensor product results from \cite{Shen.MaxInjective} and enlarges the class of positive answers to Popa's question.

The paper concludes with Section \ref{Remarks}, which contains some final remarks.  We observe that the additional ingredients Popa uses with asymptotic orthogonality to show a generator masa is maximal injective in \cite{Popa.MaxInjective} can be deduced from maximal injectivity and Ozawa's solidity of the free group factors introduced in \cite{Ozawa.Solid}. In this way, we see that these properties of the generator masas are also satisfied by the radial masa. We end with a brief discussion of maximal nuclearity of the generator and radial masas in the reduced group C$^*$-algebras of free groups.

\paragraph*{Acknowledgements}
The work in this paper originated during a visit of the fourth named author to the University of New Hampshire in April 2008. He would like to thank the faculty and students of the Department of Mathematics, and in particular Don Hadwin, for their hospitality during this visit.   Section \ref{MaxInjective} was undertaken at the Workshop in Analysis and Probability at Texas A\&{}M University in August 2008. It is a pleasure to record our gratitude to the workshop organisers and NSF for providing financial support to the workshop. 

The authors would like to thank Simon Wassermann for simulating conversations regarding maximal nuclear C$^*$-subalgebras of C$^*$-algebras.  Finally, the authors would like to thank the referee for their helpful suggestions, which have improved the exposition of the paper.

\section{Asymptotic Orthogonality and Maximal Injectivity}\label{MaxInjective}
In this section we examine how asymptotic orthogonality properties give rise to maximal injective von Neumann algebras.
\begin{definition}\label{AOP}
Let $A$ be a type $\mathrm{I}$ von Neumann subalgebra of a type \IIi von Neumann algebra $M$ with a fixed faithful normal trace $\tau_M$. Let $N$ be a finite von Neumann algebra with a fixed faithful normal trace $\tau_N$. Say that $A\subset M$ has the \emph{asymptotic orthogonality property after tensoring by $N$} if there is a non-principal ultrafilter $\omega\in\beta\mathbb N\setminus\mathbb N$ such that $x^{(1)}y_1\perp y_2x^{(2)}$ in $L^2((M\ \vnotimes\ N)^\omega,(\tau_M\otimes\tau_N)_\omega)$, whenever $x^{(1)},x^{(2)}$ are elements of $(A\ \vnotimes\ \mathbb C1)'\cap (M\ \vnotimes\ N)^\omega$ with $E_{(A\vnotimes N)^\omega}(x^{(i)})=0$ for $i=1,2$, and $y_1,y_2\in M\ \vnotimes\ N$ with $E_{A\vnotimes N}(y_i)=0$ for $i=1,2$.  Say that $A$ has the \emph{asymptotic orthogonality property} if it has this property when $N=\mathbb C1$.
\end{definition}
\noindent Note that the ultraproduct $(M\ \vnotimes\ N)^\omega$ used in this definition is constructed with respect to the product trace $\tau_M\otimes\tau_N$.

In section 2 of \cite{Popa.MaxInjective}, Popa shows that the generator masa inside a free group factor has the asymptotic orthogonality property.  The arguments given in \cite{Popa.MaxInjective} (see also \cite[Section 14.2]{Sinclair.MasaBook}) also show that the generator masa has the asymptotic orthogonality property after tensoring by any $N$ without further work. As noted in \cite{Popa.MaxInjective}, given any separable diffuse type $\mathrm{I}$ von Neumann algebra $B\neq\mathbb C1$, one can form the free product $B*R$ (where $R$ is the hyperfinite II$_1$ factor). The methods of \cite{Popa.MaxInjective} show that $B$ is maximal injective in $B*R$. Indeed, these inclusions $B\subset B*R$ have the asymptotic orthogonality property after tensoring by any $N$.

In this section we show how the maximal injectivity of a singular masa satisfying the asymptotic orthogonality property can be established without further hypotheses (Corollary \ref{MaxInjective.Masa}).  We also examine tensor products, showing how the asymptotic orthogonality property after tensoring by $N$ can be used to show that certain tensor products $A\vnotimes B$ of maximal injective subalgebras are again maximal injective. In subsequent sections we shall show that the radial masa in a free group factor also has the asymptotic orthogonality property after tensoring by any $N$. It will follow that finite tensor products of generator and radial masas are maximal injective.  We need some technical observations, the first of which easily follows from Popa's intertwining lemma.

\begin{lemma}\label{MaxInjective.Intertwine}
Let $L$ be an injective type \IIi von Neumann algebra with a separable predual equipped with a fixed faithful normal trace. Let $A$ be a type $\mathrm{I}$ von Neumann subalgebra of $L$. Then there exists a unitary $u\in L'\cap L^\omega$ with $E_{A^\omega}(u)=0$.  
\end{lemma}
\begin{proof}
Since $L$ is injective with separable predual, we can find a sequence $(L_n)$ of finite dimensional subalgebras of $L$ whose union is weakly dense in $L$.  It suffices to show that for all $\varepsilon>0$ and $n\in\mathbb N$ there exists a unitary $u\in L_n'\cap L$ with $\nm{E_A(u)}_2<\varepsilon$. However, if there exists some $n$ and $\varepsilon>0$ for which no such unitary could be found, then Popa's intertwining lemma \cite{Popa.StrongRigidity1} (see also \cite[Theorem F.12, $4\Rightarrow 1$]{BrownOzawa} for the exact statement we are using) shows that a corner of  $L_n'\cap L$ embeds into $A$ inside $L$.  This can not happen, since $L_n'\cap L$ is necessarily type \IIi and $A$ is type $\mathrm{I}_{\text{fin}}$.
\end{proof}

If $A$ is a maximal injective von Neumann subalgebra in a finite von Neumann algebra $M$, then $A$ is singular in $M$ (see \cite[Lemma 3.6]{Ozawa.OneCartan} for example). Any singular masa $A$ in $M$ necessarily has $A'\cap M\subseteq A$ so that any intermediate subalgebra $A\subset L\subset M$ has $Z(L)\subseteq L'\cap M\subseteq A'\cap M\subseteq A$.

\begin{corollary}\label{MaxInjective.Masa}
Let $A$ be a singular masa in a \IIi von Neumann algebra $M$ with a separable predual.  Suppose that $A$ has the asymptotic orthogonality property, then $A$ is maximal injective.
\end{corollary}
\begin{proof}
Let $L$ be an injective von Neumann algebra with $A\subset L\subset M$.  Let $p$ be the maximal central projection in $L$ so that $Lp$ is type $\mathrm{II}_1$. As noted above, $p\in A$. If $p\neq 0$, use Lemma \ref{MaxInjective.Intertwine} to find a unitary $u$ in $(Lp)'\cap (Lp)^\omega$ with $E_{A^\omega}(u)=0$. Choosing some non-zero $v\in Lp$ with $E_A(v)=0$, the asymptotic orthogonality property gives $uv\perp vu$. On the other hand $uv=vu$, so $uv=0$ and $v=0$, giving a contradiction. Hence $p=0$ and $L$ is a finite type $\mathrm{I}$ von Neumann algebra. By \cite{Kadison.DiagMatrices}, $A$ is regular in $L$ (see also \cite[Lemma 2.3]{Shen.MaxInjective}), so by singularity of $A$ in $M$ it follows that $L=A$.\end{proof}

The proof of our second technical observation, which allows us to handle tensor products, is contained in \cite{Shen.MaxInjective} modulo a theorem from \cite{Stratila.Commutation}. We include the details for completeness.
\begin{lemma}\label{MaxInjective.CE}
Let $A\subset M$ and $B\subset N$ be two inclusions of finite von Neumann algebras with separable preduals. Suppose that $A$ is type $\mathrm{I}$, $A'\cap M\subset A$ and that $B$ is a maximal injective von Neumann subalgebra of $N$. Let $A\ \vnotimes\ B\subset L\subset M\ \vnotimes\ N$ be an injective intermediate von Neumann subalgebra. Then
$$
E_{A\vnotimes N}(L)=A\ \vnotimes\ B.
$$
\end{lemma}
\begin{proof}
Since $A$ is finite and type $\mathrm{I}$, it follows that $(A'\cap M)'\cap M=A$. Indeed if $A$ is of type $\mathrm{I}_n$, we can express $A$ as $Z\otimes\mathbb M_n$ for some abelian $Z$ and construct a corrresponding factorisation $M=\tilde{M}\otimes\mathbb M_n$.  Then $(A'\cap M)=(Z'\cap \tilde{M})\otimes\mathbb C1$ so the condition that $A'\cap M\subset A$ implies that $Z$ is a masa in $\tilde{M}$. Hence $(A'\cap M)'\cap M=(Z\otimes \mathbb C1)'\cap (\tilde{M}\otimes\mathbb M_n)=Z\otimes\mathbb M_n=A$.  The general case follows by taking a direct sum.

Given $x\in L$, $E_{A\vnotimes N}(x)$ is the unique element of minimal $\nm{\cdot}_2$ in $\overline{\mathrm{co}}^2\{uxu^*:u\in\mathcal U((A'\cap M)\ \vnotimes\ Z(N))\}$ as $((A'\cap M)\ \vnotimes\ Z(N))'\cap (M\ \vnotimes\ N)=A\ \vnotimes\ N$.  Since $(A'\cap M)\ \vnotimes\ Z(N)\subset A\ \vnotimes\ B\subset L$, it follows that $E_{A\vnotimes N}(L)\subset L$. If $x\in L$ has $E_{A\vnotimes N}(x)\notin A\ \vnotimes\ B$, then we have
$$
A\ \vnotimes\ B\subsetneq (A\ \vnotimes\ B,E_{A\vnotimes N}(x))''\subseteq A\ \vnotimes\ N,
$$
where $(A\ \vnotimes\ B,E_{A\ \vnotimes\ N}(x))''$ is a subalgebra of $L$ and hence injective.  This contradicts \cite[Theorem 6.7]{Stratila.Commutation}, which shows that $A\ \vnotimes\ B$ is maximal injective in $A\ \vnotimes\ N$, whenever $A$ is injective and $B$ is maximal injective in a von Neumann algebra $N$ with separable predual.
\end{proof}

These ingredients enable us to gain control over certain intermediate injective subalgebras.
\begin{lemma}\label{MaxInjective.Main}
Let $M,N$ be finite von Neumann algebras with separable preduals. Let $A$ be a type $\mathrm{I}$ von Neumann subalgebra of $M$ with $A'\cap M\subseteq A$ and the asymptotic orthogonality property after tensoring by $N$. Let $B$ be a type $\mathrm{I}$ von Neumann subalgebra of $N$ which is a maximal injective von Neumann subalgebra of $N$. If $L$ is an intermediate injective von Neumann algebra between $A\ \vnotimes\ B$ and $M\ \vnotimes\ N$, then $L$ is necessarily type $\mathrm{I}$.
\end{lemma}
\begin{proof}
Let $p$ be the maximal central projection of $L$ such that $Lp$ is type $\mathrm{II}_1$. Note that $p\in L'\cap (M\ \vnotimes\ N)\subset A\ \vnotimes\ B$. If $p\neq 0$, then Lemma \ref{MaxInjective.Intertwine} gives us a unitary $u\in (Lp)'\cap (Lp)^\omega$, with $E_{((A\vnotimes B)p)^\omega}(u)=0$, where the ultrapower $(Lp)^\omega$ is constructed from the product trace on $M\ \vnotimes\ N$. Since $(A\ \vnotimes\ B)p\subsetneq Lp$, there is some non-zero $v\in Lp$ with $E_{A\vnotimes B}(v)=0$. Regarding $u\in (M\ \vnotimes\ N)^\omega$ and $v\in M\ \vnotimes\ N$, we have
$$
E_{A\vnotimes N}(v)=E_{A\vnotimes B}E_{A\vnotimes N}(v)=E_{A\vnotimes N}E_{A\vnotimes B}(v)=0,
$$
where the first equality uses Lemma \ref{MaxInjective.CE} to see that $E_{A\vnotimes N}(v)\in A\ \vnotimes\ B$.  By writing $u=(u_n)$ for some $u_n\in L$ with $E_{A\vnotimes B}(u_n)=0$ for all $n$, it follows that $E_{A\vnotimes N}(u_n)=0$ for each $n$ just as above and so $E_{(A\vnotimes N)^\omega}(u)=0$.  As $A\subset M$ has the asymptotic orthogonality property after tensoring by $N$, we have $uv\perp vu$. Since $u\in (Lp)'\cap (Lp)^\omega$ and $v\in Lp$, $u$ and $v$ commute.  Therefore $uv=0$ and hence $v=0$, which is a contradiction. Therefore $p=0$.
\end{proof}

With the notation of the previous lemma, if we additionally assume that $A$ is a singular masa with the asymptotic orthogonality property and $B$ is a masa which is maximal injective, we can then deduce the maximal injectivity of $A\ \vnotimes\ B$ in $M\ \vnotimes\ N$ as follows. Suppose $L$ is an injective von Neumann algebra between $A\ \vnotimes\ B$ and $M\ \vnotimes\ N$. By Lemma \ref{MaxInjective.Main}, $L$ is a finite type $\mathrm{I}$ von Neumann algebra.  Since $B$ is singular in $M$, $A\ \vnotimes\ B$ is singular in $M\ \vnotimes\ N$ by \cite[Corollary 2.4]{Saw.StrongSing}.  On the other hand, every masa in a finite type $\mathrm{I}$ von Neumann algebra is regular by \cite{Kadison.DiagMatrices} (see also \cite[Lemma 2.3]{Shen.MaxInjective}). Ergo $L=A\ \vnotimes\ B$.

Returning to the more general situation of type $\mathrm{I}$ algebras $A$ and $B$, the argument of the previous paragraph breaks down as there need not be a non-trivial normalising unitary of $A\ \vnotimes\ B$ in $L$. Indeed, take $P$ to be $\mathbb M_2(\mathbb C)\oplus\mathbb C1$, naturally embedded inside $Q=\mathbb M_3(\mathbb C)$, then $P'\cap Q\subseteq P$ but every normalising unitary of $P$ in $Q$ lies in $P$. To circumvent this difficulty, we need to use \emph{groupoid normalisers}. Recall that if $P\subset Q$ is an inclusion of von Neumann algebras with $P'\cap Q\subseteq P$, then a partial isometry $v\in Q$ is a \emph{groupoid normaliser} of $P$ if $vPv^*\subseteq P$ and $v^*Pv\subseteq P$. We write $\mathcal{GN}_Q(P)$ for the collection of all groupoid normalisers. Recall too that if $A$ is a singular masa in a finite von Neumann algebra $M$, then $\mathcal{GN}_M(A)\subset A$ (see \cite[Lemma 6.2.3(iv)]{Sinclair.MasaBook}). To use the argument of the preceding paragraph in a groupoid normaliser context, we need the groupoid normaliser analogue of Kadison's fact that masas in finite type $\mathrm{I}$ von Neumann algebras are regular.
\begin{lemma}\label{Type1.1}
Let $P\subset Q$ be an inclusion of finite type $\mathrm{I}$ von Neumann algebras with separable preduals and $P'\cap Q\subseteq P$. Then $\mathcal{GN}_Q(P)''=Q$.
\end{lemma}
\begin{proof}
As $P'\cap Q\subseteq P$, it follows that $Z(Q)\subseteq Z(P)$. Hence we can assume that $Q$ is type $\mathrm{I}_n$ for some $n\in\mathbb N$, as the general case follows by taking a direct sum. By \cite[Theorem 3.3]{Popa.Kadison}, $P$ contains a masa of $Q$, say $B$. By \cite[Theorem 3.19]{Kadison.DiagMatrices}, we can write $Q=A\otimes\mathbb M_n$, where $A$ is the centre of $Q$ and $B=A\otimes\mathbb D_n$, where $\mathbb D_n$ are the diagonal matrices of $\mathbb M_n$.  An easy calculation shows that $p\otimes e_{i,j}$ is a groupoid normaliser of $P$ for every projection $p\in A$ and the $(e_{i,j})_{i,j=1}^n$ are the standard matrix units of $\mathbb M_n$.  These elements evidentally generate $Q$.
\end{proof}

In \cite{Fang.CompletelySingular}, the second author introduced the notation of \emph{complete singularity} of an inclusion $P\subseteq Q$, which implies that $\mathcal{GN}_Q(P)\subset P$. Since it was shown in \cite[Proposition 3.2]{Fang.CompletelySingular} that a maximal injective von Neumann subalgebra is completely singular, the next lemma follows immediately.  Alternatively, one can establish the lemma directly from Connes' characterisation of amenability in \cite{Connes.InjectiveFactors} by showing that if $P$ is an injective von Neumann algebra and $v$ is a groupoid normaliser of $P$, then $(P\cup\{v\})''$ is also injective.
\begin{lemma}\label{GN}
Let $P$ be a maximal injective von Neumann subalgebra of a von Neumann algebra $Q$. Then $\mathcal{GN}_Q(P)\subset P$.  
\end{lemma}

We are now in a position to prove the main result of this section.
\begin{theorem}\label{MaxInjective.Final}
Let $M,N$ be finite von Neumann algebras with separable preduals. Let $A$ be a type $\mathrm{I}$ von Neumann subalgebra of $M$ with $\mathcal{GN}_M(A)\subset A$ and the asymptotic orthogonality property after tensoring by $N$. Let $B$ be a type $\mathrm{I}$ von Neumann subalgebra of $N$ which is a maximal injective von Neumann subalgebra of $N$.  Then $A\ \vnotimes\ B$ is maximal injective in $M\ \vnotimes\ N$.
\end{theorem}
\begin{proof}
Let $L$ be an injective von Neumann algebra with $A\ \vnotimes\ B\subsetneq L\subseteq M\ \vnotimes\ N$. By Lemma \ref{MaxInjective.Main}, $L$ is type $\mathrm{I}$ and so $A\ \vnotimes\ B\subset L$ is an inclusion of finite type $\mathrm{I}$ von Neumann algebras with $(A\ \vnotimes\ B)'\cap L\subset A\ \vnotimes\ B$.   By Lemma \ref{Type1.1}, there is a groupoid normaliser $v$ of $A\ \vnotimes\ B$ in $L$ with $v\notin A\ \vnotimes\ B$. On the other hand $\mathcal{GN}_M(A)\subset A$ by hypothesis and $\mathcal{GN}_N(B)\subset B$ by Lemma \ref{GN}. Corollary 5.6 of \cite{Saw.Groupoid} shows that 
$$
\mathcal{GN}_{M\vnotimes N}(A\ \vnotimes\ B)''=\mathcal{GN}_M(A)''\ \vnotimes\ \mathcal{GN}_{N}(B)''=A\ \vnotimes\ B,
$$
and this gives a contradiction. Hence $A\ \vnotimes\ B$ is maximal injective in $M\ \vnotimes\ N$.
\end{proof}

Taking $B=N=\mathbb C$ in the previous theorem immediately yields the next corollary, which we could also establish directly in exactly the same way as Corollary \ref{MaxInjective.Masa}.
\begin{corollary}
Let $A\subset M$ be an inclusion of a type $\mathrm{I}$ von Neumann algebra inside a finite von Neumann algebra with separable predual such that $\mathcal{GN}_M(A)\subset A$ and $A$ has the asymptotic orthogonality property. Then $A$ is maximal injective in $M$.
\end{corollary}

Finally we note that the asymptotic orthogonality property does not pass to tensor products.  Indeed, no inclusion $A\ \vnotimes\ B\subset M\ \vnotimes\ N$ with $A$ and $B$ masas in type \IIi factors can have the asymptotic orthogonality property.  To see this, take a unitary $u$ in $A'\cap M^\omega$ with $E_{A^\omega}(u)=0$ by Lemma \ref{MaxInjective.Intertwine} and some non-zero $v\in N$ with $E_B(v)=0$. Then $(u\otimes 1)$ lies in $(A\ \vnotimes\ B)'\cap (M\ \vnotimes N)^\omega$ with $E_{(A\vnotimes B)^\omega}(u\otimes 1)=0$ and $1\otimes v\in M\ \vnotimes\ N$ with $E_{A\vnotimes\ B}(1\otimes v)=0$ and $uv=vu$.

\section{The radial masa in a free group factor}\label{Prelim}
For $K\geq 2$, let $\F_K$ denote the free group on $K$ generators. Write $a_1\dots,a_K$ for the generators of $\F_K$.  We regard $\F_K$ as a subset of the free group factor $\LF_K$. The unique faithful normal tracial state $\tau$ on $\LF_K$ induces a pre-Hilbert space norm $\nm{x}_2=\tau(x^*x)^{1/2}$ on $\LF_K$ and completing $\LF_K$ in this norm yields the Hilbert space $\ell^2(\F_K)$. Therefore, we can regard $\LF_K$ as a subset of $\ell^2(\F_K)$ as well as an von Neumann algebra acting on $\ell^2(\F_K)$.

Given $g\in\F_K$, write $|g|$ for the length of $g$. For $n\geq 0$, let $w_n\in\LF_K$ be the sum of all words of length $n$ in $\F_K$. The recurrence relations
\begin{equation}\label{Prelim.Rec}
w_1w_n=w_nw_1=\begin{cases}w_{n+1}+(2K-1)w_{n-1},&n>1;\\w_2+2Kw_0,&n=1;\\w_1,&n=0;\end{cases}
\end{equation}
from \cite[Theorem 1]{Cohen.NormFreeGroup}, show that $w_1$ generates an abelian von Neumann subalgebra of $\LF_K$, which we denote by $A$. This is the \emph{radial} or \emph{Laplacian} subalgebra of $\LF_K$ which was shown to be masa in $\LF_K$ by Pytlik, \cite{Pytlik.RadialMasa}.  Write $L^2(A)$ for the closure of $A$ in $\ell^2(\F_K)$ and note that $(w_n/\nm{w_n}_2)_{n=0}^\infty$ forms an orthonormal basis for $L^2(A)$.

As promised in the introduction, we note that the radial and generator masas are not inner conjugate in $\LF_K$.
\begin{proposition}\label{Inner}
There is no unitary $u\in\LF_K$ with $uAu^*=\{a_1\}''$.
\end{proposition}
\begin{proof}
Let $\Phi$ be the automorphism of $\LF_K$ induced by the automorphism of $\F_K$ which interchanges the generators $a_1$ and $a_2$ and fixes the other generators. Note that $\Phi(w_1)=w_1$ so that $\Phi$ fixes the radial masa pointwise. If such a unitary $u$ existed, then $\Phi(u)A\Phi(u)^*=\Phi(uAu^*)=\Phi(\{a_1\}'')=\{a_2\}''$ and then $\theta=\mathrm{ad}\ u\Phi(u)^*$ gives an inner automorphism of $\LF_K$ with $\theta(\{a_2\}'')=\{a_1\}''$ and this contradicts \cite[Corollary 4.3]{Popa.Orth}.
\end{proof}

To show that the radial masa is singular in $\LF_K$, R\u{a}dulescu decomposed $\ell^2(\F_2)$ into a direct sum of orthogonal $A-A$-bimodules, \cite{Radulescu.LapMasa}.  Before setting out this decomposition, we need some preliminary notation. For $l\geq 0$, let $W_l$ denote the subspace of $\ell^2(\F_K)$ spanned by the words of length $l$ in $\F_K$ and let $q_l$ be the orthogonal projection from $\ell^2(\F_K)$ onto $W_l$.  Given a vector $\xi\in W_l$ for $l\geq 1$, and integers $r,s\geq 0$, define
\begin{equation}\label{Prelim.1}
\xi_{r,s} =\frac{q_{r+s+l}(w_r \xi w_s)}{(2K-1)^{(r+s)/2}}.
\end{equation}
R\u{a}dulescu's definition does not have the scaling factor $(2K-1)^{(r+s)/2}$ on the denominator, which we introduce simplify our subsequent calculations. Given a word $g$ of length $l\geq 1$, $q_{r+s+l}(w_rgw_s)$ is precisely the sum of all reduced words of length $r+s+l$ which contain $g$ from the $(r+1)$-th letter to the $(r+l)$-th letter. Since there are precisely $(2K-1)^{(r+s)}$ such words, the chosen scaling factor ensures that $\|g_{r,s}\|_2=\|g\|_2$. This observation also shows that $g_{r,s}\perp h_{r,s}$ for two distinct words $g,h$ of length $l$. As such our scaling factor ensures that $\nm{\xi_{r,s}}_2=\nm{\xi}_2$, for all $l\geq 1$, $\xi\in W_l$ and $r,s\geq 0$. Thus $\xi\mapsto\xi_{r,s}$ extends to an isometry $\ell^2(\F_K)\ominus L^2(A)\rightarrow\ell^2(\F_K)\ominus L^2(A)$. Define $\xi_{r,s}=0$ when either $r<0$ or $s<0$.  In \cite{Radulescu.LapMasa}, R\u{a}dulescu constructed a basis for $\ell^2(\F_K)\ominus L^2(A)$ which lends itself to calculations involving the radial masa and as such plays a key role in this paper.  The next lemma sets out the properties we need. Since our normalisation conventions differ from \cite{Radulescu.LapMasa}, we indicate how to derive these properties.
\begin{lemma}\label{NewLem}
There is a sequence of orthonormal vectors $(\xi^i)_{i=1}^\infty$ in $\ell^2(\F_K)$ with the following properties.
\begin{enumerate}
\item Each $\xi^i$ lies in $W_{l(i)}$ for some $l(i)\geq 1$ and has $(\xi^i)^*=\pm\xi^i$.
\item The subspaces $\overline{\mathrm{Span }(A \xi^i A)}^{2}$ are pairwise orthogonal in $\ell^2(\F_K)$ and $\oplus\overline{\mathrm{Span }(A \xi^i A)}^{2}=\ell^2(\F_K)\ominus L^2(A)$.
\item For those $i$ with $l(i)>1$, the spaces $\overline{\mathrm{Span }(A \xi^i A)}^{2}$ have orthonormal bases $(\xi^i_{r,s})_{r,s\geq 0}$.
\item For each $i,j>0$, the map $T_{i,j}$ which sends $\xi^i_{r,s}$ to $\xi^j_{r,s}$ extends to a bounded invertible operator from $\overline{\mathrm{Span }(A \xi^i A)}^{2}$ onto $\overline{\mathrm{Span }(A \xi^j A)}^{2}$. Furthermore, there exists a constant $C_0>1$ (which depends on $K$ but not $i$ or $j$) such that 
\begin{equation}\label{Prelim.2}
\nm{T_{i,j}},\nm{T_{i,j}^{-1}}\leq C_0,\quad i,j>0.
\end{equation}
\end{enumerate}
\end{lemma}
\begin{proof}
For each $l\geq 1$, let $S_l$ denote the linear subspace of $W_l$ spanned by $\{q_l(w_1g),\ q_l(gw_1):g\in\F_K,\ |g|\leq l-1\}$.  As in Theorem 7 of \cite{Radulescu.LapMasa}, for each $l\geq 1$ we choose an orthonormal basis $(\xi^{j,l})_j$ for $W_l\ominus S_l$. As the subspaces $W_l\ominus S_l$ are self-adjoint, we can additionally insist that the elements of this basis satisfy $(\xi^{j,l})^*=\pm\xi^{j,l}$ (as R\u{a}dulescu does when $l=1$ in \cite[Theorem 7]{Radulescu.LapMasa}). Re-indexing these elements as $(\xi^i)_i$, gives a sequence $(\xi^i)_i$, where $\xi^i\in W_{l(i)}$ which satisfies condition (i).  Condition (ii) follows from \cite[Lemma 4]{Radulescu.LapMasa} (the first statement is part (b) of this lemma, while part (a) gives the second statement).  Lemma 3(a) of \cite{Radulescu.LapMasa} ensures that for those $i$ with $l(i)\geq 2$, the elements $(\xi^{i}_{r,s})_{r,s\geq 0}$ are pairwise orthogonal. Since our normalisation conventions above ensure that $\|\xi^{i}_{r,s}\|_2=\|\xi^{i}\|_2$ for all such $r,s$ condition (iii) follows. 

When $l(i),l(j)>1$ the maps $T_{i,j}$ in condition (iv) are unitary operators by (iii). It then suffices to establish (iv) for some constant $C_0'$ when $l(i)=1$ and some fixed $j$ with $l(j)>1$ as all other cases follow by composition.  In this case our map $T_{i,j}$ is the map $T_0'$ appearing in the proof of \cite[Lemma 6]{Radulescu.LapMasa} and which is noted to be bounded and invertible there. One can explicitly estimate the bound on these maps from the results in \cite{Radulescu.LapMasa}, however this is not necessary as there are only finitely many $i$ with $l(i)=1$, so there is certainly a uniform bound $C_0'$ on the norm of $T_{i,j}$ and $T_{i,j}^{-1}$.
\end{proof}

Note that the set $\{\xi^i_{r,s}:i\geq 1,r,s\geq 0\}$ does not give an orthonormal basis for $\ell^2(\F_K)\ominus L^2(A)$, due to the presence of some indices $i$ with $l(i)=1$. However, the last two facts above show that this set is at least a basis for $\ell^2(\F_K)\ominus L^2(A)$ and the $2$-norm it induces is equivalent to the norm $\nm{\cdot}_2$.  We shall subsequently use this in a tensor product setting.
\begin{proposition}\label{Prelim.Estimate} 
Let $\mathcal H$ be a Hilbert space and let $\eta\in(\ell^2(\F_K)\otimes\mathcal H)\ominus(L^2(A)\otimes\mathcal H)$. Then there exist vectors $(\alpha_{r,s}^i)_{r,s\geq 0,i>0}$ in $\mathcal H$ such that ${\displaystyle \eta=
\sum_{\substack{i\geq 1\\r\geq 0,s\geq 0}}  \xi^{i}_{r,s}\otimes\alpha^{i}_{r,s}}.$ Furthermore
\[C_0^{-1}\nm{\eta}_{\ell^2(\F_K)\otimes \mathcal H}\leq \Big(\sum_{\substack{i\geq 1\\r\geq 0,s\geq 0}}
\|\alpha^{i}_{r,s}\|_{\mathcal H}^2\Big)^{1/2}\leq C_0\nm{\eta}_{\ell^2(\F_K)\otimes\mathcal H},
\]
where $C_0>1$ is the constant appearing in equation (\ref{Prelim.2}).
\end{proposition}
\begin{remark}\label{NewRem}
When the Hilbert space $\mathcal H$ is of the form $L^2(N)$ for some finite von Neumann algebra $N$, it makes sense to refer to the conjugate $x^*$ of an element $x\in \ell^2(\F_K)\otimes L^2(N)$ (indeed, this is defined by extending the conjugation operation from $\LF_K\ \vnotimes\ N$ to the Hilbert space). We can expand both $x$ and $x^*$ in terms of the $\xi^i_{r,s}$, writing $x=\sum_{i\geq 1,\ r,s\geq 0}\xi^i_{r,s}\otimes\alpha^i_{r,s}$ and $x^*=\sum_{i\geq 1,\ r,s\geq 0}\xi^i_{r,s}\otimes\beta^i_{r,s}$. Since Lemma \ref{NewLem} (i) and (\ref{Prelim.1}) give $\xi^i_{r,s}=\pm\xi^i_{s,r}$, it follows that $\beta^i_{r,s}=\pm\alpha^i_{s,r}$ for each $i,r$ and $s$.
\end{remark}

A critical component of R\u{a}dulescu's computations in \cite{Radulescu.LapMasa} are recurrence relations analogous to (\ref{Prelim.Rec}) showing how the $\xi^i_{r,s}$ behave under multiplication with $w_1$. We will need these relations subsequently so we recall them here. Again it is slightly more convenient for our purposes to use a suitable normalisation, so we define $\tilde{w_1}=w_1/(2K-1)^{1/2}$. Applying our normalisation conventions to \cite[Lemma 1]{Radulescu.LapMasa}, for $i\geq 1$ and $r,s\geq 0$, there exist values of $\sigma(i)=\pm 1$ so that
\begin{equation}\label{Prelim.Rec.2}
\tilde{w_1}\xi^i_{r,s}=\begin{cases}\xi^i_{r+1,s}+\xi^i_{r-1,s},&r\geq 1;\\ \xi^i_{1,s},&r=0,\ l(i)\geq 2;\\\xi^i_{1,s}+\frac{\sigma(i)}{2K-1}\xi_{0,s-1}^i,&r=0,\ l(i)=1;\end{cases}
\end{equation}
and
\begin{equation}\label{Prelim.Rec.3}
 \xi^i_{r,s}\tilde{w_1}=\begin{cases}\xi^i_{r,s+1}+\xi^i_{r,s-1},&s\geq 1;\\\xi^i_{r,1},&s=0,\ l(i)\geq 2;\\\xi^i_{r,1}+\frac{\sigma(i)}{2K-1}\xi^i_{r-1,0},&s=0,\ l(i)=1.\end{cases}
\end{equation}
Since $\xi^i_{-1,s}$ and $\xi^i_{r,-1}$ are defined to be zero, the second cases of the relations above are identical to the first cases.

\section{Locating the support of elements of $A'\cap(\LF_K)^\omega$}\label{Asymptotics}
In this section we establish the first half of the technical estimates required to show that the radial masa has the asymptotic orthogonality property in $\LF_K$.  Our objective, realised in Lemma \ref{Asymptotics.Main}, is to control the support of elements of $A'\cap(\LF_K)^\omega$ which are orthogonal to $A^\omega$. We find that the essential support of elements in $\LF_K$ orthogonal to $A$ and approxiately commuting with $A$ must be contained in the closed linear span of the $\xi^i_{r,s}$ for sufficiently large $r$ and $s$.  We perform these calculations after tensoring by an arbitrary finite von Neumann algebra $N$ as our objective is the asymptotic orthogonality property after tensoring by such an $N$. For each $i$, define
$$
\lambda(i)=\begin{cases}1,&l(i)\geq 2;\\1-\frac{\sigma(i)}{2K-1},&l(i)=1.\end{cases}
$$
\begin{lemma}\label{Asymptotics.Tech1}
Let $N$ be a finite von Neumann algebra and suppose that $x\in\ell^2(\F_K)\otimes L^2(N)$ has $x\perp L^2(A\ \vnotimes\ N)$ and $\nm{x}_2=1$. Write 
$$
{\displaystyle x=\sum_{r,s\geq 0,\ i\geq 1}\xi^i_{r,s}\otimes\alpha^{i}_{r,s}},
$$
with convergence in $\ell^2(\F_K)\otimes L^2(N)$, for some $\alpha^i_{r,s}\in L^2(N)$. Then, for each $s'\geq s\geq 1$ 
\begin{align}
&\Big|\Big(\sum_{\substack{r\geq s'\\i\geq 1}}\|\alpha_{r-s,0}^{i}+\lambda(i)\alpha_{r-s+2,0}^{i}+\dots+\lambda(i)\alpha_{r+s-2,0}^{i}+\lambda(i)\alpha_{r+s,0}^{i}\|_2^2\Big)^{1/2}-\Big(\sum_{\substack{r\geq s'\\i\geq 1}}\|\alpha_{r,s}^{i}\|_2^2\Big)^{1/2}\Big|\nonumber\\\
&\leq\Big(\sum_{\substack{r\geq s\\i\geq 1}}\|\alpha^i_{r,s}-(\alpha_{r-s,0}^{i}+\lambda(i)\alpha_{r-s+2,0}^{i}+\dots+\lambda(i)\alpha_{r+s-2,0}^{i}+\lambda(i)\alpha_{r+s,0}^{i})\|_2^2\Big)^{1/2}\nonumber\\
&\leq 3^{s-1}C_0\nm{x(\tilde{w_1}\otimes 1)-(\tilde{w_1}\otimes 1)x}_2,\label{Asymptotics.Support.4}
\end{align}
where $C_0$ is the constant of Proposition \ref{Prelim.Estimate}.
\end{lemma}
\begin{proof}
Note that the first inequality in the Lemma is immediate from the triangle inequality, so it suffices to prove the second inequality. Write $\varepsilon=\nm{x(\tilde{w_1}\otimes 1)-(\tilde{w_1}\otimes 1)x}_2$ and define $\alpha^i_{r,s}=0$ whenever $r<0$ or $s<0$.  Recalling the convention that $\xi^i_{r,s}$ is zero whenever $r<0$ or $s<0$, the recurrence relations (\ref{Prelim.Rec.2}) and (\ref{Prelim.Rec.3}) show that
$$
[\tilde{w_1},\xi^i_{r,s}]=\begin{cases}\xi^i_{r+1,s}+\xi^i_{r-1,s}-\xi^i_{r,s+1}-\xi^i_{r,s-1},& l(i)\geq 2,\ r,s\geq 0\\\xi^i_{r+1,s}+\xi^i_{r-1,s}-\xi^i_{r,s+1}-\xi^i_{r,s-1},&l(i)=1,\ r,s\geq 1
\\\xi^i_{1,s}-\xi^i_{0,s+1}-\lambda(i)\xi^i_{0,s-1},&l(i)=1,\ r=0,\ s>0\\\xi^i_{r+1,0}+\lambda(i)\xi^i_{r-1,0}-\xi^i_{r,1},&l(i)=1,\ r>0,\ s=0\\
\xi^i_{1,0}-\xi^i_{0,1},&l(i)=1,\ r=s=0.
\end{cases}
$$
These relations give
\begin{align}
[\tilde{w_1}\otimes 1,x]=&\sum_{r,s\geq 0,\ l(i)\geq 2}(\xi^r_{r+1,s}+\xi^i_{r-1,s}-\xi^i_{r,s+1}-\xi^i_{r,s-1})\otimes\alpha_{r,s}^i\label{neweq}\\
&+\sum_{r,s\geq 1,\ l(i)=1}(\xi^i_{r+1,s}+\xi^i_{r-1,s}-\xi^i_{r,s+1}-\xi^i_{r,s-1})\otimes \alpha^i_{r,s}\nonumber\\
&+\sum_{s\geq 1,\ l(i)=1}(\xi^i_{1,s}-\xi^i_{0,s+1}-\lambda(i)\xi^i_{0,s-1})\otimes\alpha_{0,s}^i\nonumber\\
&+\sum_{r\geq 1,\ l(i)1=1}(\xi^i_{r+1,0}+\lambda(i)\xi^i_{r-1,0}-\xi^i_{r,1})\otimes\alpha_{r,0}^i\nonumber\\
&+\sum_{l(i)=1}(\xi_{1,0}^i-\xi^i_{0,1})\otimes\alpha_{0,0}^i\nonumber,
\end{align}
with convergence in $\ell^2(\F_K)\otimes L^2(N)$.  We then rearrange this sum to write it in the form $\sum_{r,s\geq 0,\ i\geq 1}\xi^i_{r,s}\otimes\beta^i_{r,s}$ for some $\beta^i_{r,s}\in L^2(N)$.  For $l(i)\geq 2$ and $r,s\geq 0$ contributions to $\beta^i_{r,s}$ arise from the $(r+1,s)$, $(r-1,s)$, $(r,s+1)$ and $(r,s-1)$ terms in the first sum on the right of (\ref{neweq}) (the cases when either $r=0$ or $s=0$ do not cause extra difficulty, as our notational conventions ensure that $\xi^i_{-1,s}=\xi^i_{r,-1}=0$ and $\alpha^i_{-1,s}=\alpha^i_{r,-1}=0$). Thus $\beta_{r,s}^i=\alpha_{r-1,s}^{i}+\alpha_{r+1,s}^{i}-\alpha_{r,s-1}^{i}-\alpha_{r,s+1}^{i}$ when $l(i)=2$.  The values of $\beta^i_{r,s}$ when $l(i)=1$ can also be computed from (\ref{neweq}). For example we can only obtain a contribution to $\beta^i_{0,0}$ with $l(i)=1$ from the $s=1$ term in line 3 of (\ref{neweq}) and the $r=1$ term of line 4, giving $\beta^i_{0,0}=\lambda(i)(\alpha^i_{1,0}-\alpha^i_{0,1})$ in this case.  Continuing in this fashion, we obtain
\begin{align*}
[\tilde{w_1}\otimes 1,x]=&\sum_{\substack{l(i)\geq 2,\ r,s\geq 0}}\xi^i_{r,s}\otimes(\alpha_{r-1,s}^{i}+\alpha_{r+1,s}^{i}-\alpha_{r,s-1}^{i}-\alpha_{r,s+1}^{i})\\
&+\sum_{\substack{l(i)=1,\ r,s\geq 1}}\xi^i_{r,s}\otimes(\alpha_{r-1,s}^{i}+\alpha_{r+1,s}^{i}-\alpha_{r,s-1}^{i}-\alpha_{r,s+1}^{i})\\
&+\sum_{\substack{l(i)=1,\ r\geq 1}}\xi^i_{r,0}\otimes (\lambda(i)\alpha^i_{r+1,0}+\alpha^i_{r-1,0}-\alpha^i_{r,1})\\
&+\sum_{l(i)=1,\ s\geq 1}\xi^i_{0,s}\otimes(\alpha^i_{1,s}-\alpha_{0,s-1}^i-\lambda(i)\alpha_{0,s+1}^i)\\
&+\sum_{l(i)=1}\lambda(i)\xi^i_{0,0}\otimes(\alpha^i_{1,0}-\alpha^i_{0,1}).
\end{align*}
Proposition \ref{Prelim.Estimate} gives
\begin{align}
\label{Asymptotics.Support.Rad1}C_0^2\varepsilon^2\geq&\sum_{\substack{l(i)\geq 2,\ r,s\geq 0}}\|\alpha_{r-1,s}^{i}+\alpha_{r+1,s}^{i}-\alpha_{r,s-1}^{i}-\alpha_{r,s+1}^{i}\|_2^2\\
&+\sum_{\substack{l(i)=1,\ r,s\geq 1}}\|\alpha_{r-1,s}^{i}+\alpha_{r+1,s}^{i}-\alpha_{r,s-1}^{i}-\alpha_{r,s+1}^{i}\|_2^2\nonumber\\
&+\sum_{\substack{l(i)=1,\ r\geq 1}}\|\lambda(i)\alpha^i_{r+1,0}+\alpha^i_{r-1,0}-\alpha^i_{r,1}\|_2^2\nonumber\\
&+\sum_{l(i)=1,\ s\geq 1}\|\alpha^i_{1,s}-\alpha_{0,s-1}^i-\lambda(i)\alpha_{0,s+1}^i\|_2^2\nonumber\\
&+\sum_{l(i)=1}\|\lambda(i)(\alpha^i_{1,0}-\alpha^i_{0,1})\|_2^2\nonumber,
\end{align}
where $C_0$ is the constant appearing in Proposition \ref{Prelim.Estimate}.

Considering the terms with $s=0$ from the first line of right hand side of (\ref{Asymptotics.Support.Rad1}) and those terms from the third line gives
$$
\sum_{\substack{r\geq 1\\l(i)\geq 2}}\|\alpha_{r-1,0}^{i}+\alpha_{r+1,0}^{i}-\alpha_{r,1}^{i}\|_2^2+\sum_{\substack{r\geq 1\\l(i)=1}}\|\lambda(i)\alpha^i_{r+1,0}+\alpha^i_{r-1,0}-\alpha^i_{r,1}\|_2^2\leq C_0^2\varepsilon^2,
$$
since $\alpha_{r,-1}^{i}$ is defined to be zero. As $\lambda(i)=1$ whenever $l(i)\geq 2$, we have 
$$
\sum_{\substack{r\geq 1\\i\geq 1}}\|\alpha^i_{r-1,0}+\lambda(i)\alpha^i_{r+1,0}-\alpha^i_{r,1}\|_2^2\leq C_0^2\varepsilon^2.
$$
This gives the $s=1$ case of (\ref{Asymptotics.Support.4}).

For the $s=2$ case, take those terms on the right-hand side of (\ref{Asymptotics.Support.Rad1}) with $s=1$ and $r\geq 2$ to obtain
\begin{equation}\label{Asymptotics.Support.7}
\sum_{\substack{r\geq 2\\i\geq 1}}\|\alpha_{r-1,1}^{i}+\alpha_{r+1,1}^{i}-\alpha_{r,0}^{i}-\alpha_{r,2}^{i}\|_2^2\leq C_0^2\varepsilon^2.
\end{equation}
Now
\begin{align}
\alpha_{r,2}^{i}+\alpha_{r,0}^{i}-\alpha_{r+1,1}^{i}-\alpha_{r-1,1}^{i}=&\Big(\alpha_{r,2}^{i}-(\alpha_{r-2,0}^{i}+\lambda(i)\alpha_{r,0}^{i}+\lambda(i)\alpha_{r+2,0}^{i})\Big)\label{Asymptotic.Support.8}\\&-\Big(\alpha_{r+1,1}^{i}-(\alpha_{r,0}^{i}+\lambda(i)\alpha_{r+2,0}^{i})\Big)\nonumber\\&-\Big(\alpha_{r-1,1}^{i}-(\alpha_{r-2,0}^{i}+\lambda(i)\alpha_{r,0}^{i})\Big),\nonumber
\end{align}
so, (\ref{Asymptotics.Support.7}) and the $s=1$ case of the lemma combine to give
\begin{align}
&\Big(\sum_{r\geq 2,\ i\geq 1}\|\alpha_{r,2}^{i}-(\alpha_{r-2,0}^{i}+\lambda(i)\alpha_{r,0}^{i}+\lambda(i)\alpha_{r+2,0}^{i})\|_2^2\Big)^{1/2}\nonumber\\
\leq&C_0\varepsilon+\Big(\sum_{r\geq 2,\ i\geq 1}\|\alpha^i_{r+1,1}-(\alpha^i_{r,0}+\lambda(i)\alpha^i_{r+2,0})\|_2^2\Big)^{1/2}+\Big(\sum_{r\geq 2,\ i\geq 1}\|\alpha_{r-1,1}^{i}-(\alpha_{r-2,0}^{i}+\lambda(i)\alpha_{r,0}^{i}\|_2^2\Big)^{1/2}\label{NewEq.1}\\
\leq&C_0\varepsilon+C_0\varepsilon+C_0\varepsilon=3C_0\varepsilon\nonumber.
\end{align}
We need to sum over $r\geq 2$ above so that the sum in the second term in (\ref{NewEq.1}) is over those $r+1\geq 1$ which allows the $s=1$ case of the lemma to be used.

Suppose inductively that the lemma holds for $1\leq s\leq s_0$ with $s_0\geq 2$. Taking those terms in (\ref{Asymptotics.Support.Rad1}) with $s=s_0$ and $r\geq s_0$ gives
\begin{equation}\label{Asymptotics.Support.5}
\sum_{\substack{r\geq s_0\\i\geq 1}}\left\|\alpha_{r-1,s_0}^{i}+\alpha_{r+1,s_0}^{i}-\alpha_{r,s_0-1}^{i}-\alpha_{r,s_0+1}^{i}\right\|_2^2\leq C_0^2\varepsilon^2.
\end{equation}
The identity
\begin{align}
&\alpha_{r,s_0-1}^{i}+\alpha_{r,s_0+1}^{i}-\alpha_{r-1,s_0}^{i}-\alpha_{r+1,s_0}^{i}\label{Asymptotics.Support.6}\\
=&\alpha_{r,s_0+1}^{i}-\left(\alpha_{r-(s_0+1),0}^{i}+\lambda(i)\alpha_{r-(s_0-2),0}^{i}+\dots+\lambda(i)\alpha_{r+s_0-2,0}^{i}+\lambda(i)\alpha_{r+(s_0+1),0}^{i}\right)\nonumber\\
&+\alpha_{r,s_0-1}^{i}-\left(\alpha_{r-(s_0-1),0}^{i}+\lambda(i)\alpha_{r-(s_0-4),0}^{i}+\dots+\lambda(i)\alpha_{r+s_0-4,0}^{i}+\lambda(i)\alpha_{r+(s_0-1),0}^{i}\right)\nonumber\\
&-\left(\alpha_{r+1,s_0}^{i}-\left(\alpha_{r+1-s_0),0}^{i}+\lambda(i)\alpha_{r+1-(s_0-3),0}^{i}+\dots+\lambda(i)\alpha_{r+1+(s_0-3),0}^{i}+\lambda(i)\alpha_{r+1+s_0),0}^{i}\right)\right)\nonumber\\
&-\left(\alpha_{r-1,s_0}^{i}-\left(\alpha_{r-1-s_0,0}^{i}+\lambda(i)\alpha_{r-1-(s_0-3),0}^{i}+\dots+\lambda(i)\alpha_{r-1+(s_0-3),0}^{i}+\lambda(i)\alpha_{r-1+s_0,0}^{i}\right)\right)\nonumber,
\end{align}
the inductive hypothesis, (\ref{Asymptotics.Support.5}) and the triangle inequality imply that
\begin{align*}
&\Big(\sum_{\substack{r\geq s_0+1\\i\geq 1}}\|\alpha_{r,s_0+1}^{i}-(\alpha_{r-(s_0+1),0}^{i}+\lambda(i)\alpha_{r-(s_0-2),0}^{i}+\dots+\lambda(i)\alpha_{r+s_0-2,0}^{i}+\lambda(i)\alpha_{r+(s_0+1),0}^{i})\|_2^2\Big)^{1/2}\\
\leq&C_0\varepsilon+3^{(s_0-2)}C_0\varepsilon+3^{(s_0-1)}C_0\varepsilon+3^{(s_0-1)}C_0^2\varepsilon\leq  3\cdot 3^{(s_0-1)}C_0\varepsilon=3^{s_0}C_0\varepsilon.
\end{align*}
By induction, the lemma holds for all $s$.
\end{proof}

The crux of Lemma \ref{Asymptotics.Main} is contained in the next Lemma.
\begin{lemma}\label{Asymptotics.Tech2}
Let $N$ be a finite von Neumann algebra. Let $x=(x_n)$ be an element of $(A\otimes\mathbb C1)'\cap (\LF_K\otimes N)^\omega$ with $\nm{x_n}_2=1$ and $E_{A\vnotimes N}(x_n)=0$ for all $n$. Write $x_n=\sum_{i,r,s}\xi^i_{r,s}\otimes \alpha_{r,s}^{n,i}$ for some $\alpha_{r,s}^{n,i}\in L^2(N)$ and convergence in $\ell^2(\F_K)\otimes L^2(N)$. Then
$$
\lim_{n\rightarrow\omega}\sum_{i\geq 1,\ r\geq 0}\|\alpha^{n,i}_{r,0}\|_2^2=\lim_{n\rightarrow\omega}\sum_{i\geq 1,\ s\geq 0}\|\alpha^{n,i}_{0,s}\|_2^2=0.
$$
\end{lemma}
\begin{proof}
Let $\varepsilon_n=\nm{x_n(\tilde{w_1}\otimes 1)-(\tilde{w_1}\otimes 1)x_n}_2$.  For $t\geq 1$,
\begin{align}
&\big(\sum_{r\geq 0,\ i\geq 1}\|\alpha^{n,i}_{r,0}+\lambda(i)\alpha_{r+2,0}^{n,i}+\dots+\lambda(i)\alpha_{r+2t,0}^{n,i}\|_2^2\big)^{1/2}\nonumber\\
=&\big(\sum_{r+t\geq t,\ i\geq 1}\|\alpha^{n,i}_{(r+t)-t,0}+\lambda(i)\alpha_{r+t-(t-2),0}^{n,i}+\dots+\lambda(i)\alpha_{(r+t)+t,0}^{n,i}\|_2^2\big)^{1/2}\nonumber\\
\leq&\big(\sum_{r+t\geq t,\ i\geq 1}\|\alpha_{r+t,t}^{n,i}\|_2^2\big)^{1/2}+3^{t-1}C_0\varepsilon_n=\big(\sum_{r\geq t,\ i\geq 1}\|\alpha_{r,t}^{n,i}\|_2^2\big)^{1/2}+3^{t-1}C_0\varepsilon_n,\label{NewEq.2}
\end{align}
where the inequality comes from Lemma \ref{Asymptotics.Tech1}. Another use of this lemma with $s'=t+1$ and $s=t-1$ gives
\begin{align}
&\big(\sum_{r\geq 0,\ i\geq 1}\|\alpha^{n,i}_{r+2,0}+\lambda(i)\alpha_{r+4,0}^{n,i}+\dots+\lambda(i)\alpha_{r+2t,0}^{n,i}\|_2^2\big)^{1/2}\nonumber\\
=&\big(\sum_{r+t+1\geq t+1,\ i\geq 1}\|\alpha^{n,i}_{(r+t+1)-(t-1),0}+\lambda(i)\alpha_{(r+t+1)-(t-3),0}^{n,i}+\dots+\lambda(i)\alpha_{(r+t+1)+(t-1),0}^{n,i}\|_2^2\big)^{1/2}\nonumber\\
\leq&\big(\sum_{r+t+1\geq t+1,\ i\geq 1}\|\alpha_{r+t+1,t-1}^{n,i}\|_2^2\big)^{1/2}+3^{t-2}C_0\varepsilon_n=\big(\sum_{r\geq t+1,\ i\geq 1}\|\alpha_{r,t-1}^{n,i}\|_2^2\big)^{1/2}+3^{t-2}C_0\varepsilon_n\label{NewEq.3}.
\end{align}
Using
$$
\frac{1}{2}\|\beta-\gamma\|_2^2\leq\|\beta\|_2^2+\|\gamma\|_2^2,\quad \beta,\gamma\in L^2(N),
$$
with $\beta=\alpha^{n,i}_{r,0}+\lambda(i)\alpha^{n,i}_{r+2,0}+\dots+\lambda(i)\alpha^{n,i}_{r+2t,0}$ and $\gamma=\alpha^{n,i}_{r+2,0}+\lambda(i)\alpha^{n,i}_{r+4,0}+\dots+\lambda(i)\alpha^{n,i}_{r+2t,0}$, the inequalities (\ref{NewEq.2}) and (\ref{NewEq.3}) give
\begin{align}
&\frac{1}{2}\sum_{\substack{r\geq 0,\ i\geq 1}}\|\alpha^{n,i}_{r,0}-(1-\lambda(i))\alpha^{n,i}_{r+2,0}\|_2^2\label{Asymptotics.Tech2.20}\\
\leq&\Big(\big(\sum_{r\geq t,\ i\geq 1}\|\alpha_{r,t}^{n,i}\|_2^2\big)^{1/2}+3^{t-1}C_0\varepsilon_n\Big)^2+\Big(\big(\sum_{r\geq t+1,\ i\geq 1}\|\alpha_{r,t-1}^{n,i}\|_2^2\big)^{1/2}+3^{t-2}C_0\varepsilon_n\Big)^2\nonumber\\
\leq&\sum_{\substack{r\geq t\\i\geq 1}}\|\alpha^{n,i}_{r,t}\|_2^2+\sum_{\substack{r\geq t+1\\i\geq 1}}\|\alpha^{n,i}_{r,t-1}\|_2^2\nonumber\\
&+\Big(2\cdot 3^{t-1}C_0(\sum_{\substack{r\geq t\\i\geq 1}}\|\alpha^{n,i}_{r,t}\|_2^2)^{1/2}+2\cdot 3^{t-2}C_0(\sum_{\substack{r\geq t+1\\i\geq 1}}\|\alpha^{n,i}_{r,t-1}\|_2^2)^{1/2}\Big)\varepsilon_n+(3^{2t-2}+3^{2t-4})C_0^2\varepsilon_n^2\nonumber\\
=&\sum_{\substack{r\geq t\\i\geq 1}}\|\alpha^{n,i}_{r,t}\|_2^2+\sum_{\substack{r\geq t+1\\i\geq 1}}\|\alpha^{n,i}_{r,t-1}\|_2^2\nonumber\\
&+C_0^2\Big(2\cdot 3^{t-1}+2\cdot 3^{t-2}\Big)\varepsilon_n+C_0^2(3^{2t-2}+3^{2t-4})\varepsilon_n^2,\nonumber
\end{align}
for each $n,t\in\mathbb N$, where the estimate  $\sum_{r,s,i}\|\alpha^{n,i}_{r,s}\|_2^2\leq C_0^2\|x_n\|_2^2=C_0^2$ comes from Proposition \ref{Prelim.Estimate}.

Summing (\ref{Asymptotics.Tech2.20}) over $t=1,\dots,t_0$, and crudely estimating the first two sums on the right-hand side gives
\begin{align*}
&\frac{t_0}{2}\sum_{\substack{r\geq 0,\ i\geq 1}}\|\alpha^{n,i}_{r,0}-(1-\lambda(i))\alpha^{n,i}_{r+2,0}\|_2^2\\
&\leq2\sum_{\substack{r\geq 0,t\geq 0\\i\geq 1}}\|\alpha^{n,i}_{r,t}\|_2^2+C_0^2\sum_{t=1}^{t_0}(2\varepsilon_n(3^{t-1}+3^{t-2})+(3^{2t-2}+3^{2t-4})\varepsilon_n^2)
\end{align*}
Using $\sum_{r,t,i}\|\alpha_{r,t}^{n,i}\|_2^2\leq C_0^2\|x_n\|_2^2=C_0^2$ again, we have\begin{align}
&\sum_{\substack{r\geq 0,\ i\geq 1}}\|\alpha^{n,i}_{r,0}-(1-\lambda(i))\alpha^{n,i}_{r+2,0}\|_2^2\nonumber\\
&\leq\frac{4C_0^2}{t_0}+\frac{2C_0^2}{t_0}\sum_{t=1}^{t_0}\left(2\varepsilon_n(3^{t-1}+3^{t-2})+(3^{2t-2}+3^{2t-4})\varepsilon_n^2\right).\label{Asymptotics.Tech2.1}
\end{align}
The first term of (\ref{Asymptotics.Tech2.1}) can be made arbitrarily small by fixing a suitably large value for $t_0$. Thus
\begin{equation}\label{Asymptotics.Tech2.10}
\lim_{n\rightarrow\omega}\sum_{\substack{r\geq 0,\ i\geq 1}}\|\alpha^{n,i}_{r,0}-(1-\lambda(i))\alpha^{n,i}_{r+2,0}\|^2=0,
\end{equation}
as $\lim_{n\rightarrow\omega}\varepsilon_n=0$.

When $l(i)>1$, we have $\lambda(i)=1$ so (\ref{Asymptotics.Tech2.10}) implies
$$
\lim_{n\rightarrow\omega}\sum_{\substack{r\geq 0,\ l(i)>1}}\|\alpha^{n,i}_{r,0}\|_2^2=0.
$$
For those $i$ with $l(i)=1$, we have $|(1-\lambda(i))|=1/(2K-1)$.  Since
\begin{align*}
&\left(1-\frac{1}{2K-1}\right)\|\alpha_{r,0}^{n,i}\|_2^2-\frac{1}{2K-1}\|\alpha_{r+2,0}^{n,i}\|_2^2\\
&\leq\|\alpha_{r,0}^{n,i}\|_2^2-\frac{2}{2K-1}|\langle\alpha_{r,0}^{n,i},\alpha_{r+2,0}^{n,i}\rangle|\\
&\leq\|\alpha_{r,0}^{n,i}+(1-\lambda(i))\alpha^{n,i}_{r+2,0}\|_2^2,
\end{align*}
(\ref{Asymptotics.Tech2.10}) also implies that
$$
\left(1-\frac{2}{2K-1}\right)\sum_{r\geq 0,\ l(i)=1}\|\alpha^{n,i}_{r,0}\|_2^2\leq\sum_{r\geq 0,\ l(i)=1}\|\alpha_{r,0}^{n,i}-(1-\lambda(i))\alpha^{n,i}_{r+2,0}\|_2^2\rightarrow 0,\quad n\rightarrow\omega.
$$
Combining the cases $l(i)=1$ and $l(i)>1$ gives
$$
\lim_{n\rightarrow\omega}\sum_{\substack{r\geq 0,\ i\geq 1}}\|\alpha^{n,i}_{r,0}\|_2^2=0.
$$

Finally, replacing $x$ by $x^*$ interchanges the roles of $r$ and $s$ (see Remark \ref{NewRem}). Therefore
$$
\lim_{n\rightarrow\omega}\sum_{\substack{s\geq 0,\ i\geq 1}}\|\alpha^{n,i}_{0,s}\|_2^2=0,
$$
which is the second limit required for the lemma.
\end{proof}

We can now deduce the main result of this section, showing that the elements of $\LF_K$ which are orthogonal to $A$ and approximately commute with $w_1$ are essentially supported on the span of the $\xi^i_{r,s}$ for large $r$ and $s$.  
\begin{lemma}\label{Asymptotics.Main}
Let $N$ be a finite von Neumann algebra. Let $x=(x_n)$ be an element of $(A\otimes\mathbb C1)'\cap (\LF_K\otimes N)^\omega$ with $\nm{x_n}_2=1$ and $E_{A\vnotimes N}(x_n)=0$ for all $n$. Write $x_n=\sum_{i,r,s}\xi^i_{r,s}\otimes \alpha_{r,s}^{n,i}$ for some $\alpha_{r,s}^{n,i}\in L^2(N)$ and convergence in $\ell^2(\F_K)\otimes L^2(N)$.  For each $m\in\mathbb N$,
$$
\lim_{n\rightarrow\omega}\sum_{\substack{i\geq 1,\\r\leq m\text{ or }s\leq m}}\|\alpha^{n,i}_{r,s}\|_2^2=0.
$$
\end{lemma}
\begin{proof}
Lemma \ref{Asymptotics.Tech2} gives
\begin{equation}\label{Asymptotics.Main.2}
\lim_{n\rightarrow\omega}\sum_{\substack{i\geq 1\\r\geq 0}}\|\alpha^{n,i}_{r,0}\|_2^2=0.
\end{equation}
For fixed $s_0>0$, we can combine (\ref{Asymptotics.Main.2}) with Lemma \ref{Asymptotics.Tech1} and the triangle inequality to obtain
\begin{equation}\label{Asymptotics.Main.1}
\lim_{n\rightarrow\omega}\sum_{\substack{i\geq 1\\r\geq s_0}}\|\alpha^{n,i}_{r,s_0}\|_2^2=0.
\end{equation}
By replacing $x$ with $x^*$ (using Remark \ref{NewRem}) we interchange the roles of $r$ and $s$ and hence, for each $r_0\geq 0$
\begin{equation}\label{Asymptotics.Main.3}
\lim_{n\rightarrow\omega}\sum_{\substack{i\geq 1\\s\geq r_0}}\|\alpha^{n,i}_{r_0,s}\|_2^2=0.
\end{equation}
Combining the limits (\ref{Asymptotics.Main.1}) and (\ref{Asymptotics.Main.3}) establishes the lemma.
\end{proof}

\section{Counting Words in $\F_K$}\label{Count}
This section contains the combinatorial ingredients required to show that the radial masa has the asymptotic orthogonality property.

Let $g,h\in\F_K$. Say that there are exactly $i\geq 0$ \emph{cancelations in the product $gh$} if $|gh|=|g|+|h|-2i$. Given non-empty subsets $\sigma,\tau$ of $\{a_1^{\pm 1},\dots,a_K^{\pm 1}\}$ and $n\geq 0$, let $w_n(\sigma,\tau)$ denote the sum of all words in $\F_K$ of length $n$ which begin with an element of $\sigma$ and end with an element of $\tau$.  Write $\nu_n(\sigma,\tau)$ for the number of words in this sum so that $\nu_n(\sigma,\tau)=\nm{w_n(\sigma,\tau)}_2^2$. We abuse notation when sets of singletons are involved, writing $\nu_n(a_1,\tau)$ instead of $\nu_n(\{a_1\},\tau)$, for example. The values $\nu_n(\sigma,\tau)$ can be explicitly computed by solving certain difference equations, see Lemma 4.3 and Remark 4.4 of \cite{Sinclair.LapMasa},  but for our purposes the following estimate of Sinclair and Smith will suffice.
\begin{proposition}[({\cite[Corollary 4.5]{Sinclair.LapMasa}})]\label{Count.SS}
There exists a constant $C_1>0$, which depends only on $K$, such that 
$$
\left|\nu_n(\sigma_1,\tau_1)-\nu_n(\sigma_2,\tau_2)\right|\leq C_1,
$$
for all $n\geq 0$ and all subsets $\sigma_1,\sigma_2,\tau_1,\tau_2\subset\{a_1^{\pm 1},\dots,a_K^{\pm 1}\}$ with $|\sigma_1|=|\sigma_2|$ and $|\tau_1|=|\tau_2|$.
\end{proposition}

Recall that for a vector $\xi\in W_l$, and $m\geq 0$ we defined $(\xi)_{m,m}=(2K-1)^{-m}q_{2m+l}(w_m\xi w_m)$ in (\ref{Prelim.1}). In this section, our objective is to estimate
$$
\Big|\ip{(k_1)_{m,m}(g_1-g_2)}{h(k_2)_{m,m}}\Big|,
$$
where $g_1,g_2,h,k_1,k_2\in\F_K$ have $|g_1|=|g_2|$, $k_1\neq e$, $k_2\neq e$ and $m$ is sufficiently large. Fix such $g_1,g_2,h,k_1$ and $k_2$ and let $m>2\max(|g_1|,|h|)$.

Write $S^m(k_1,g_1)$ for the collection of all words of the form $xk_1yg_1$, where $x,y$ have length $m$ and there are no cancelations in the products $xk_1$ and $k_1y$. Write $T^m(k_2,h)$ for the collection of all words of the form $hxk_2y$, where $x,y$ have length $m$ and there are no cancelations in the products $xk_2$ and $k_2y$.  These definitions are constructed so that
$$
(k_1)_{m,m}g_1=\frac{1}{(2K-1)^m}\sum_{s\in S^m(k_1,g_1)}s,\quad h(k_2)_{m,m}=\frac{1}{(2K-1)^m}\sum_{t\in T^m(k_2,h)}t,
$$
as $(k_i)_{m,m}$ is the sum of all words $xk_iy$ with $|x|=|y|=m$ and no cancelations in the products $xk_i$ and $k_iy$ normalised by a factor of $(2K-1)^{-m}$.
Therefore
\begin{equation}\label{Count.2}
\ip{(k_1)_{m,m}g_1}{h(k_2)_{m,m}}=\frac{1}{(2K-1)^{2m}}|T^m(k_2,h)\cap S^m(k_1,g_1)|.
\end{equation}
For $0\leq i\leq|g_1|$, write $S_i^m(k_1,g_1)$ for those words $xk_1yg_1$ in $S^m(k_1,g_1)$ which have exactly $i$ cancelations in the product $yg_1$ so that $S^m(k_1,g_1)=\bigcup_{i=0}^{|g_1|}S_i^m(k_1,g_1)$.

\begin{lemma}\label{Count.Tech.2}
With the notation above, 
$$
\Big||T^m(k_2,h)\cap S^m(k_1,g_1)|-|T^m(k_2,h)\cap S^m(k_1,g_2)|\Big|\leq C_1(2K-1)^{m+|h|},
$$
where $C_1$ is the constant of Proposition \ref{Count.SS}.
\end{lemma}
\begin{proof}
Let $l=|g_1|=|g_2|$ and express $g_1$ and $g_2$ as reduced words $g_1=g_{1,1}\dots g_{1,l}$, $g_2=g_{2,1}\dots g_{2,l}$. There are $(2K-1)^m$ reduced words $x\in\F_K$ of length $m$ with no cancelations in the product $xk_2$.  Given such an $x$, write 
$$
I_p(x)=\Big|\{y\in\F_K:hxk_2y\in S^m(k_1,g_p),\ |y|=m,\ |k_2y|=m+|k_2|\}\Big|,\quad p=1,2
$$
and it then suffices to show that
\begin{equation}\label{Count.1}
|I_1(x)-I_2(x)|\leq C_1(2K-1)^{|h|},
\end{equation}
for all such $x$. Fix such an $x$, and let $j$ be the number of cancelations in the product $hx$. By comparing the lengths of words, $I_1(x)=I_2(x)=0$ unless there is some integer $i$ with $0\leq i\leq |g_1|$ and 
$$
|k_1|+|g_1|-2i=|k_2|+|h|-2j,
$$
for some $j$ with $0\leq j\leq |h|$. Furthermore, the words $hxk_2y$ lying in the sets $S^m(k_1,g_1)$ and $S^m(k_1,g_2)$ appearing in (\ref{Count.1}) must actually lie in $S^m_i(k_1,g_1)$ and $S^m_i(k_1,g_2)$ respectively. We now consider three cases individually.

Firstly suppose that $|hxk_2|\geq m+|k_1|$. The words of $S^m(k_1,g_1)$ and $S^m(k_1,g_2)$ all contain a copy of $k_1$ from the $(m+1)$-th letter to the $(m+|k_1|)$-th letter. Therefore we can assume that $k_1$ is contained in $hxk_2$ from the $(m+1)$-th letter to the $(m+|k_1|)$-th letter, otherwise both $I_1(x)=0$ and $I_2(x)=0$ in which case (\ref{Count.1}) is immediate.  For $hxk_2y$ to lie in $S^m_i(k_1,g_1)$, the last $l-i$ letters of $y$ must be $g_{1,i+1}\dots g_{1,l}$.  That is, writing $y=y_1\dots y_m$, we must have $y_{m-l+i+1}=g_{1,i+1}, y_{m-l+i+2}=g_{1,i+2},\dots,y_m=g_{1,l}$. Hence $hxk_2y$ must have the form
$$
hxk_2y_1\dots y_{m-l+i}g_{1,i+1}\dots g_{1,l},
$$
where there is no cancelation in the product $k_2y_1$, $y_{m-l+i}\neq y_{m-l+i+1}^{-1}= g_{1,i+1}^{-1}$ (this condition gives no restriction if $i=l$) and $y_{m-l+i}\neq g_{1,i}$ (this condition gives no restriction if $i=0$). This last condition is necessary as we must have
\begin{equation}\label{New.1011}
hxk_2y_1\dots y_{m-l+i}(g_{1,i}^{-1}\dots g_{1,1}^{-1}g_{1,1}\dots g_{1,i})g_{1,i+1}\dots g_{1,l},
\end{equation}
with no cancelation between $y_{m-l+i}$ and $g_{1,i}^{-1}$ as otherwise $y_{m-l+i}$ would cancel with $g_{1,i+1}$ in (\ref{New.1011}). Furthermore every such $y_1\dots y_{m-l+i}$ gives rise to some $y$ with $hxk_2y\in S^m(k_1,g_1)$. Therefore
\begin{align*}
I_1(x)=\nu_{m-l+i}(\{a_1^{\pm 1},\dots,a_K^{\pm 1}\}\setminus\{\text{last letter of }k_2\},\{a_1^{\pm 1},\dots,a_K^{\pm 1}\}\setminus\{g_{1,i+1}^{-1},g_{1,i}\})
\end{align*}
with the appropriate adjustments if $i=0$ or $i=l$. Similarly 
$$
I_2(x)=\nu_{m-l+i}(\{a_1^{\pm 1},\dots,a_K^{\pm 1}\}\setminus\{\text{last letter of }k_2\},\{a_1^{\pm 1},\dots,a_K^{\pm 1}\}\setminus\{g_{2,i+1}^{-1},g_{2,i}\}),
$$
with the appropriate adjustments if $i=0$ or $i=l$.  Proposition \ref{Count.SS} implies that$$
|I_1(x)-I_2(x)|\leq C_1,
$$
and so (\ref{Count.1}) follows.

Now consider the case that $m+1\leq |hxk_2|<m+|k_1|$. Let $r=|hxk_2|-m$. We can assume that the first $r$ letters of $k_1$ make up the last $r$ letters of $hxk_2$, otherwise both $I_1(x)$ and $I_2(x)$ are zero.  Any $y$ for which $hxk_2y\in S_i^m(k_1,g_1)$ must commence with the last $(|k_1|-r)$-letters of $k_1$, that is we may assume $y_1,\dots,y_{|k_1|-r}$ are given. Arguing exactly as in the previous paragraph we see that (making adjustments if $i=0$ or $i=l$)
\begin{align*}
I_1(x)=\nu_{m-l-(r+1)+i}(\{a_1^{\pm 1},\dots,a_K^{\pm 1}\}\setminus\{y_{|k_1|-r}^{-1}\},\{a_1^{\pm 1},\dots,a_K^{\pm 1}\}\setminus\{g_{1,i+1}^{-1},g_{1,i}\})
\end{align*}
and
\begin{align*}
I_2(x)=\nu_{m-l-(r+1)+i}(\{a_1^{\pm 1},\dots,a_K^{\pm 1}\}\setminus\{y_{|k_1|-r}^{-1}\},\{a_1^{\pm 1},\dots,a_K^{\pm 1}\}\setminus\{g_{2,i+1}^{-1},g_{2,i}\}).
\end{align*}
Again (\ref{Count.1}) follows from Proposition \ref{Count.SS}.

Finally, consider the case that $|hxk_2|\leq m$ and let $r=m-|hxk_2|$ and note that $r= m-(m+|k_2|+|h|-2j)\leq |h|$ (recalling that $0\leq j\leq |h|$). In this case any $y$ for which $hxk_2y$ lies in $S_i^m(k_1,g_1)$ must contain a copy of $k_1$ from the $(r+1)$-th position to the $(r+|k_1|)$-th position.  Any choices of $y_1,\dots,y_r$ which lead to no cancelations in $k_2y$ and $y_rk_1$ are permissible and for each such choice, arguing just as before (and making appropriate adjustments when $i=0$ or $i=l$), there are
$$
\nu_{m-l-(r+|k_1|)+i}(\{a_1^{\pm 1},\dots,a_N^{\pm 1}\}\setminus\{y_{r+|k_1|}^{-1}\},\{a_1^{\pm 1},\dots,a_N^{\pm 1}\}\setminus\{g_{1,i+1}^{-1},g_{1,i}\})
$$
choices of $y_{r+|k_1|+1}\dots y_m$ giving rise to an element $hxk_1y$ in $S_i^m(k_1,g_1)$. There are at most $(2K-1)^{r}$ choices of $y_1,\dots,y_r$, so Proposition \ref{Count.SS} gives 
$$
|I_1(x)-I_2(x)|\leq C_1(2K-1)^r\leq C_1(2K-1)^{|h|}.
$$
Thus (\ref{Count.1}) also holds in this case.
\end{proof}

We need one final combinatorial ingredient for our proof of the maximal injectivity of the radial masa.
\begin{lemma}\label{Count.Tech.1}
Let $g_1,g_2,h\in\F_K$ have $|g_1|=|g_2|$. Write $I$ for the set of all pairs $(k_1,k_2)$ in $\F_K$ with $k_1\neq e$ and $k_2\neq e$ such that $\ip{(k_1)_{m,m}(g_1-g_2)}{h(k_2)_{m,m}}\neq 0$ for some $m>2\max(|g_1|,|h|)$. Then there exists a constant $C_2$, depending only on $K$, $|g_1|$ and $|h|$ so that, for $k_1$ fixed, there are at most $C_2$ values of $k_2$ with $(k_1,k_2)\in I$ and for $k_2$ fixed, there are at most $C_2$ values of $k_1$ with $(k_1,k_2)\in I$.
\end{lemma}
\begin{proof}
Fix $k_1\in\F_K$. For $m>2\max(|g_1|,|h|)$, the inner product $\ip{(k_1)_{m,m}(g_1-g_2)}{h(k_2)_{m,m}}$ is zero when both the intersections $S^m(k_1,g_1)\cap T^m(k_2,h)$ and $S^m(k_1,g_2)\cap T^m(k_2,h)$ are empty.  These intersections are both empty unless there are $i$ and $j$ with $0\leq i\leq |g_1|$ and $0\leq j\leq |h|$ and a word of length $2m+|k_1|+|g_1|-2i=2m+|k_2|+|h|-2j$ containing $k_1$ from the $(m+1)$-letter to the $(m+|k_1|$)-th letter and containing $k_2$ from the $(m+|h|-2j+1)$-th letter to the $(m+|k_2|+|h|-2j)$-th letter. Since the possible values of $i$ and $j$ depend only on $|g_1|$ and $|h|$, it suffices to find $C_2'$, depending only on $|g_1|$, $h$ and $K$ (and not $m$), such that for each $i$ and $j$, there are at most $C_2'$ values of $k_2$ satisfying these conditions.  For fixed $i$ and $j$, we have the freedom to choose the first $\max(0,2j-|h|)$-letters of $k_2$ and the last $\max(0,|g_1|-2i)$-letters of $k_2$ (subject to their being no cancelations) the remaining letters are determined by $k_1$. Estimating crudely, there are at most $(2K)^{|h|+|g_1|}$ such choices and so we can take $C_2'=(2K)^{|h|+|g_1|}$. A similar argument shows that the value of $C_2$ obtained also satisfies the second statement of the lemma.
\end{proof}

\section{Maximal Injectivity of the Radial Masa}\label{Radial}
In this section, we combine the results of the previous two sections to show that the radial masa has the asymptotic orthogonality property.  The next lemma converts the combinatorial arguments of the previous section into a suitable form for use in Theorem \ref{Radial.AOP}.
\begin{lemma}\label{FinalEstimate}
Let $g_1,g_2,h\in\F_K$ have $|g_1|=|g_2|$. There exists a constant $C_3$, which depends only on $g_1,g_2,h$ and $K$, such that for all finite von Neumann algebras $N$, all $m>2\max(|g_1|,|h|)$ and all vectors $\eta_1,\eta_2\in\ell^2(\F_K)\otimes L^2(N)$ which lie in the closed linear span of $\xi^i_{r+m,s+m}\otimes L^2(N)$ for $i\geq 1$ and $r,s\geq 1$, we have
$$
\Big|\ip{\eta_1((g_1-g_2)\otimes z_{1,2})}{(h\otimes z_{2,2})\eta_2}\Big|\leq C_3(2K-1)^{-m}\nm{\eta_1}_2\nm{\eta_2}_2\nm{z_{1,2}}\nm{z_{2,2}},
$$
for all $z_{1,2},z_{2,2}\in N$.
\end{lemma}
\begin{proof}
Note that $\xi^i_{r+m,s+m}=(\xi^i_{r,s})_{m,m}$. Each $\xi^i_{r,s}$ lies in $W_{l(i)+r+s}$, the span of the words of length $l(i)+r+s$ and the map $\zeta\mapsto\zeta_{m,m}$ is an isometry $W_{l(i)+r+s}\rightarrow W_{l(i)+r+s+2m}$. Thus $\xi^i_{r+m,s+m}$ lies in the span of the elements $k_{m,m}$ for $k\in\F_r$ with $|k|=l(i)+r+s$.  Note too that the elements $(k_{m,m})_{k\in\F_r\setminus\{e\}}$ form an orthonormal set in $\ell^2(\F_K)$. This follows as $k_{m,m}$ consists of the normalised sum of all words of length $|k|+2m$ which contain $k$ beginning in the $(m+1)$-th position. All the words in the sum for $k_1$ are then orthogonal to the words in the sum for $k_2$ unless $k_1=k_2$.

The hypotheses of the lemma and the preceding paragraph allow us to write
$$
\eta_1=\sum_{k\neq e}k_{m,m}\otimes\alpha_k,\quad\eta_2=\sum_{k\neq e}k_{m,m}\otimes\beta_k,
$$
for some $\alpha_k,\beta_k\in L^2(N)$ with $\sum_k\|\alpha_k\|_2^2=\nm{\eta_1}_2^2$ and $\sum_k\|\beta_k\|_2^2=\nm{\eta_2}_2^2$ and so
\begin{align}
&\Big|\ip{\eta_1((g_1-g_2)\otimes z_{1,2})}{(h\otimes z_{2,2})\eta_2}\Big|\nonumber\\
=&\left|\sum_{k_1,k_2\neq e} \langle (k_1)_{m,m}(g_1-g_2), h(k_2)_{m,m}\rangle \langle \alpha_{k_1}z_{1,2}, z_{2,2}\beta_{k_2}\rangle \right|\nonumber\\
\leq&\nm{z_{1,2}}\nm{z_{2,2}}\sum_{k_1,k_2\neq e}\|\alpha_{k_1}\|_2\|\beta_{k_2}\|_2\Big|\ip{(k_1)_{m,m}(g_1-g_2)}{h(k_2)_{m,m}}\Big|.\label{MaxInjective.Tech.1}
\end{align}
Let $I$ be the set of all pairs $(k_1,k_2)$ such that the inner product $\ip{(k_1)_{m,m}(g_1-g_2)}{h(k_2)_{m,m}}$ is non-zero for some $m>2\max(|g_1|,|h|)$. By Lemma \ref{Count.Tech.2} and (\ref{Count.2})
$$
|\ip{(k_1)_{m,m}(g_1-g_2)}{h(k_2)_{m,m}}|\leq C_1(2K-1)^{|h|-m},\quad m>2\max(|g_1|,|h|),
$$
where $C_1$ is the constant (depending only on $K$) from Proposition \ref{Count.SS}.  Applying the Cauchy-Schwartz inequality to (\ref{MaxInjective.Tech.1})) gives
\begin{align*}
&\Big|\ip{\eta_1((g_1-g_2)\otimes z_{1,2})}{(h\otimes z_{2,2})\eta_2}\Big|\\
\leq&C_1(2K-1)^{|h|-m}\big(\sum_{(k_1,k_2)\in I}\|\alpha_{k_1}\|_2^2\big)^{1/2}\big(\sum_{(k_1,k_2)\in I}\|\beta_{k_2}\|_2^2\big)^{1/2}\nm{z_{1,2}}\nm{z_{2,2}}.
\end{align*}
Lemma \ref{Count.Tech.1} gives us a constant $C_2$ (depending only on $|g_1|,|h|$ and $K$) such that for each $k_1$, there are at most $C_2$ words of $k_2$ with $(k_1,k_2)\in I$ and for each $k_2$, there are at most $C_2$ words $k_1$ with $(k_1,k_2)\in I$. Therefore
$$
\sum_{(k_1,k_2)\in I}\|\alpha_{k_1}\|_2^2\leq C_2\sum_{k_1}\|\alpha_{k_1}\|_2^2=C_2\|\eta_1\|_2^2,\text{ and } \sum_{(k_1,k_2)\in I}\|\beta_{k_2}\|_2^2\leq C_2\|\eta_2\|_2^2.
$$
Thus
$$
\Big|\ip{\eta_1((g_1-g_2)\otimes z_{1,2})}{(h\otimes z_{2,2})\eta_2}\Big|\leq C_1C_2(2K-1)^{|h|-m}\nm{\eta_1}_2\nm{\eta_2}_2\nm{z_{1,2}}\nm{z_{2,2}},
$$
and the result follows, where the constant $C_3$ is given by $C_3=C_1C_2(2K-1)^{|h|}$.
\end{proof}

We are now in a position to prove that the radial masa in the free group factor $\LF_K$ has the asymptotic orthogonality property (after tensoring by any finite von Neumann algebra).  The maximal injectivity of the radial masa then follows from the results of Section \ref{MaxInjective}.  In fact we are able to show slightly more; in the theorem which follows we do not require that $E_{A\vnotimes N}(y_2)=0$.  This is to be expected, as a similar result for the generator masa can be obtained from Popa's calculations in \cite{Popa.MaxInjective}.
\begin{theorem}\label{Radial.AOP}
Let $N$ be a finite von Neumann algebra with a fixed faithful normal, normalised trace. Let $x^{(1)},x^{(2)}\in (A\otimes\mathbb C1)'\cap (\LF_K\ \vnotimes\ N)^\omega$ and $y_1,y_2\in \LF_K\ \vnotimes\ N$ satisfy $E_{(A\vnotimes N)^\omega}(x^{(1)})=E_{(A\vnotimes N)^\omega}(x^{(2)})=0$ and $y_1,y_2\in \LF_K\ \vnotimes\ N$ with $E_{A\vnotimes N}(y_1)=0$. Then $x^{(1)}y_1\perp y_2x^{(2)}$. In particular, the radial masa $A$ in the free group factor $\LF_K$ has the asymptotic orthogonality property after tensoring by $N$.
\end{theorem}
\begin{proof}
By linearity and density, we may first assume that $y_2$ is of the form $h\otimes z_{2,2}$ for a single group element $h\in\F_K$ and $z_{2,2}\in N$.  We may also assume that $y_1$ is an elementary tensor $z_{1,1}\otimes z_{1,2}$, where $z_{1,1}$ is an element of $\mathbb C\F_K$ with $E_A(z_{1,1})=0$, and this space decomposes as $\bigoplus_{l\geq 1}W_l\ominus\mathbb Cw_l$.  Each $W_l\ominus \mathbb Cw_l$ is spanned by elements of the form $g_1-g_2$ with $|g_1|=|g_2|=1$.  To see this take some non-zero $r=\sum_{|g|=l}\beta_gg\in W_l\ominus \mathbb Cw_l$ and note that $\sum_{|g|=l}\beta_g=0$ as $r\perp w_l$. The support of $r$ consists of those $g$ with $\beta_g\neq 0$.  If this support has precisely two elements then $r$ is already a multiple of some $g_1-g_2$.  Otherwise, choose some $g_1,g_2$ in this support and define $r'=r-\beta_{g_1}(g_1-g_2)$.  This still lies in $W_l\ominus\mathbb Cw_l$ and has a strictly smaller support. Proceeding by induction, it follows that the elements $g_1-g_2$ with $|g_1|=|g_2|=l$ span $W_l\ominus \mathbb Cw_l$.

Thus we may assume that $y_1=(g_1-g_2)\otimes z_{1,2}$ for some $|g_1|=|g_2|$. Choose representatives $(x^{(1)}_n),(x^{(2)}_n)$ of $x^{(1)},x^{(2)}$ with $E_{A\vnotimes N}(x^{(1)}_n)=E_{A\vnotimes N}(x_n^{(2)})=0$ and $\|x^{(1)}_n\|_2=\|x^{(1)}\|_2$, $\|x^{(2)}_n\|_2=\|x^{(2)}\|_2$ for all $n$.  

We can write $x^{(j)}_n=\sum_{i,r,s}\xi_{r,s}^i\otimes \alpha^{j,n,i}_{r,s}$, for $\alpha^{j,n,i}_{r,s}\in L^2(N)$.  For each $m>0$, Lemma \ref{Asymptotics.Main} shows that
$$
\lim_{n\rightarrow\omega}\sum_{\substack{i\geq 1\\r\leq m\text{ or }s\leq m}}\|\alpha^{j,n,i}_{r,s}\|_2^2=0,\quad j=1,2.
$$
In particular, if we define elements 
$$
\eta^{(j)}_n=\sum_{\substack{i\geq 1\\r,s>m}}\xi^i_{r,s}\otimes\alpha_{r,s}^{j,n,i}
$$
in $\ell^2(\F_K)\otimes L^2(N)$, then
\begin{align*}
\ip{x^{(1)}((g_1-g_2)\otimes z_{1,2})}{(h\otimes z_{2,2})x^{(2)}}&=\lim_{n\rightarrow\omega}\ip{x^{(1)}_n((g_1-g_2)\otimes z_{1,2})}{(h\otimes z_{2,2})x^{(2)}_n}\\&=\lim_{n\rightarrow\omega}\ip{\eta^{(1)}_n((g_1-g_2)\otimes z_{1,2})}{(h\otimes z_{2,2})\eta^{(2)}_n}.
\end{align*}
Lemma \ref{FinalEstimate} gives us a constant $C_3$, which depends only on $|g_1|,|h|$ and $N$ so that
$$
\Big|\ip{\eta^{(1)}_n((g_1-g_2)\otimes z_{1,2})}{(h\otimes z_{2,2})\eta^{(2)}_n}\Big|\leq C_3(2K-1)^{-m}\nm{\eta^{(1)}_n}_2\nm{\eta^{(2)}_n}_2\nm{z_{1,2}}\nm{z_{2,2}}.
$$
Since $\lim_{n\rightarrow\omega}\|\eta^{(j)}_n\|_2=\|x^{(j)}\|_2$, we obtain
$$
\Big|\ip{x^{(1)}((g_1-g_2)\otimes z_{1,2})}{(h\otimes z_{2,2})x^{(2)}}\Big|\leq C_3(2K-1)^{-m}\nm{x^{(1)}}_2\nm{x^{(2)}}_2\nm{z_{1,2}}\nm{z_{2,2}},
$$
for all $m\in\mathbb N$. Hence $x^{(1)}((g_1-g_2)\otimes z_{1,2})\perp (h\otimes z_{2,2})x^{(2)}$, as required.
\end{proof}

As the radial masa is singular in $\LF_K$ by R\u{a}dulescu's original computation \cite{Radulescu.LapMasa} and so has $\mathcal{GN}_{\LF_K}(A)\subset A$ (see \cite[Lemma 6.2.3(iv)]{Sinclair.MasaBook}), the corollaries below follow from Corollary \ref{MaxInjective.Masa} and Theorem \ref{MaxInjective.Final}.
\begin{corollary}
For each $K\geq 2$, the radial masa is a maximal injective von Neumann subalgebra of $\LF_K$.
\end{corollary}

\begin{corollary}
Let $A$ be the radial masa in $\LF_K$. Let $B\subset N$ be an inclusion of a type $\mathrm{I}$ von Neumann algebra in a finite von Neumann algebra such that $B$ is maximal injective in $N$. Then $A\ \vnotimes\ B$ is maximal injective in $\LF_K\ \vnotimes\ N$.
\end{corollary}

\begin{corollary}
Fix $n\in\mathbb N$. For $1\leq i\leq n$, let $A_i$ be either a generator or radial masa in $\LF_{K_i}$ for some $K_i\geq 2$. Then the tensor product $A_1\ \vnotimes\ \dots\ \vnotimes\ A_{n}$ is maximal injective in $\LF_{K_1}\ \vnotimes\ \dots\ \vnotimes\ \LF_{K_{n}}$.
\end{corollary}

\section{Concluding Remarks}\label{Remarks}
Let $B$ denote a generator masa in $\LF_K$. As a preliminary step in his proof that $B$ is maximal injective in $\LF_K$, Popa uses the fact (see \cite[Proposition 4.1]{Popa.Orth}) that if $B_0$ is a diffuse subalgebra of $B$, then $B_0'\cap \LF_K=B$. It follows that if $L$ is any intermediate von Neumann algebra between $B$ and $\LF_K$, then there is a family of orthogonal central projections $(p_i)_{i\geq 0}$ in $L$ which sum to $1$ such that $Lp_0=Bp_0$ and for $i\geq 1$, each $Lp_i$ is a \IIi factor.  Popa then uses the critical asymptotic orthogonality calculation \cite[Lemma 2.1]{Popa.MaxInjective} to show that these factors $Lp_i$ can not have property $\Gamma$ and the maximal injectivity of $B$ follows (see Corollary 3.3 of \cite{Popa.MaxInjective}). 

The argument of section \ref{MaxInjective} allows us to deduce the maximal injectivity of $B$ directly from Popa's asymptotic orthogonality calculation. Using Ozawa's solidity of $\LF_K$ \cite{Ozawa.Solid}, the additional properties of the generator masa found in Popa's proof can then be recovered. This works for any maximal injective von Neumann subalgebra of a solid $\mathrm{II}_1$ factor. In particular, the radial masa shares all the properties of the generator masas developed in \cite{Popa.MaxInjective}. Recall that a von Neumann algebra $M$ is \emph{solid} if the relative commutant of every diffuse von Neumann subalgebra of $M$ is injective.
\begin{proposition}
Let $B$ be a maximal injective von Neumann subalgebra of a solid \IIi factor $M$. If $L$ is an intermediate von Neumann subalgebra between $B$ and $M$, then there is a family of orthogonal central projections $(p_i)_{i\geq 0}$ in $L$ which sum to $1$ such that $Lp_0=Bp_0$ and for $i\geq 1$, each $Lp_i$ is a \IIi factor without property $\Gamma$.
\end{proposition}
\begin{proof}
Note that $B'\cap M\subseteq B$ by maximal injectivity, so that $L'\cap M=Z(L)\subseteq Z(B)$. Let $p_0$ be the maximal central projection of $L$ such that $Lp_0$ is diffuse.  Then $Lp_0\subseteq (Z(L)p_0\oplus B(1-p_0))'\cap M$ and this last algebra is injective by solidity of $M$. Since $(Z(L)p_0\oplus B(1-p_0))'\cap M$ contains $B$, we have $Lp_0\subseteq B$ by maximal injectivity of $B$.  Let $(p_i)_{i>0}$ denote the minimal projections of $Z(L)(1-p_0)$ so that for $i>0$, each $Lp_i$ is a non injective factor of type $\mathrm{II}_1$.  These factors can not have property $\Gamma$, by Popa's observation \cite[Proposition 7]{Ozawa.Solid}.
\end{proof}

The generator and radial algebras also give rise to masas at the level of the reduced C$^*$-algebra.  Write $B_0$ for the C$^*$-subalgebra of $C^*_r(\F_K)$ generated by the first generator and $A_0$ for the C$^*$-subaglebra of $C^*_r(\F_K)$ generated by $w_1$. Both $A_0$ and $B_0$ are maximal abelian subalgebras of $C^*_r(\F_K)$. The later claim is well known, while the former can be found in \cite{Pytlik.RadialMasa}.
\begin{proposition}
The algebras $A_0$ and $B_0$ are maximal nuclear C$^*$-subalgebras of $C^*_r(\F_K)$.
\end{proposition}
\begin{proof}
Represent $C^*_r(\F_K)$ on $\ell^2(\F_K)$ and let $L_0$ be a nuclear C$^*$-subalgebra of $C^*_r(\F_K)$ which contains $A_0$ or $B_0$.  Then $L_0''$ is an injective von Neumann subalgebra of $\LF_K$ containing the radial masa $A$ or the generator masa $B$ so that $L_0''=A$ or $L_0''=B$ by maximal injectivity of these algebras.  In particular $L_0''$ and hence $L_0$ is abelian. Since $A_0$ and $B_0$ are maximal abelian subalgebras of $C^*_r(\F_K)$ it follows that $L_0=A_0$ or $L_0=B_0$.
\end{proof}

An identical argument shows that any finite tensor product of $n$ copies of the algebras $A_0$ or $B_0$ is maximal nuclear in the spatial tensor product of $n$ copies of $C^*_r(\F_K$), since \cite[Theorem 4]{Wassermann.SliceMap} the tensor product of masas is maximal abelian in the spatial tensor product. As the algebra $B_0$ has the extension property (\cite[Example (i)]{Archbold.ExtensionStates}), it follows that $B_0\otimes B_0$ is a masa in $C^*_r(\F_K)\otimes_\alpha C^*_r(\F_K)$ for any C$^*$-norm $\alpha$, \cite{Wassermann.TensorMasa}.  However, $B_0\otimes B_0$ need not be a maximal nuclear C$^*$-algebra in $C^*_r(\F_K)\otimes_\alpha C^*_r(\F_K)$ when $\alpha$ is larger than the minimal norm. We would like to thank Simon Wassermann bringing this to our attention and for allowing us to include the following argument here.

Let $\alpha$ denote the conjugacy norm on $C^*_r(\F_K)\odot C^*_r(\F_K)$, given by the representation $\lambda\times\rho$ of $C^*_r(\F_K)\odot C^*_r(\F_K)$ on $\mathbb B(\ell^2(\F_K))$, where $\lambda$ and $\rho$ are the left and right regular representations respectively.  The kernel of the canonical $^*$-homomorphism $\pi:C^*_r(\F_K)\otimes_\alpha C^*_r(\F_K)\rightarrow C^*_r(\F_K)\otimes C^*_r(\F_K)$ consists precisely of the compact operators $\mathbb K(\ell^2(\F_K))$ and this forms the only non-trivial ideal in this tensor product \cite{Akemann.TensorProduct} (see also \cite{Wassermann.TensorFreeGroup}). Consequently the C$^*$-subalgebra $D$ of $C^*_r(\F_K)\otimes_\alpha C^*_r(\F_K)$ generated by $B_0\otimes B_0$ and $\mathbb K(\ell^2(\F_K))$ is a nuclear extension of a nuclear C$^*$-algebra so nuclear. Furthermore, this algebra $D$ is a maximal nuclear C$^*$-subalgebra of $C^*_r(\F_K)\otimes_\alpha C^*_r(\F_K)$. Indeed, if $D_1$ is a nuclear C$^*$-subalgebra of $C^*_r(\F_K)\otimes_\alpha C^*_r(\F_K)$ containing $D$, then $\pi(D_1)$ is a nuclear C$^*$-subalgebra containing $B_0\otimes B_0$ in $C^*_r(\F_K)\otimes C^*_r(\F_K)$. By maximal nuclearity of $B_0\otimes B_0$ in $C^*_r(\F_K)\otimes C^*_r(\F_K)$, we have $\pi(D_1)=B_0\otimes B_0$. Thus $D_1=D$.

\end{document}